\documentclass[11pt]{article}
\usepackage{mathrsfs}
\usepackage{epic,eepic,epsf,epsfig}
\usepackage{amsfonts,srcltx,mathrsfs}
\usepackage{amsmath}
\usepackage{amssymb}
\usepackage{amsbsy}
\usepackage{graphicx}
\usepackage{amsfonts}
\usepackage{color}
\usepackage{mathrsfs}
\usepackage{epic,eepic,epsf,epsfig}
\usepackage{graphicx}
\usepackage{amsfonts,srcltx,mathrsfs}
\usepackage{array}
\usepackage{amsthm}
\newtheorem{theorem}{Theorem}
\newtheorem{lem}{Lemma}
\newtheorem{cor}{Corollary}
\newtheorem{defi}{Definition}

\def\f{\noindent}

\begin{document}

\markboth{ et. al}{The generalized connectivity of some regular graphs}

\title{The generalized connectivity of some regular graphs
}
\author{Shu-Li Zhao, Rong-Xia Hao\footnote{Corresponding author. Email: 17118434@bjtu.edu.cn, rxhao@bjtu.edu.cn (R.-X. Hao) }\\[0.2cm]
{\em\small Department of Mathematics, Beijing Jiaotong University,}\\ {\small\em
Beijing 100044, P.R. China}}

\date{}
\maketitle

 The generalized $k$-connectivity $\kappa_{k}(G)$ of a graph $G$ is a parameter that can measure the reliability of a network $G$ to connect any $k$ vertices in $G$, which is proved to be NP-complete for a general graph $G$. 
 In this paper, we study the generalized $3$-connectivity of some general $m$-regular and $m$-connected graphs $G_{n}$ constructed recursively and obtain that $\kappa_{3}(G_{n})=m-1$, which attains the upper bound of $\kappa_{3}(G)$ [Discrete Mathematics 310 (2010) 2147-2163] given by Li {\em et al.} for $G=G_{n}$. As applications of the main result, the generalized $3$-connectivity of many famous networks such as the alternating group graph $AG_{n}$, the $k$-ary $n$-cube $Q_{n}^{k}$, the split-star network $S_{n}^{2}$ and the bubble-sort-star graph $BS_{n}$ etc. can be obtained directly.

\medskip

\f {\em Keywords:} Interconnection network; Generalized connectivity; Fault-tolerance; Regular graph.

\section{Introduction}
For an interconnection network, one mainly concerns the reliability and fault tolerance. An interconnection network is usually modeled as a connected graph $G=(V, E)$, where nodes represent processors
and edges represent communication links between processors.  {\it The connectivity $\kappa (G)$} of a graph $G$ is defined as the minimum number of vertices whose deletion results in a disconnected graph, which is an important parameter to evaluate the reliability and fault tolerance of a network. In addition, Whitney~\cite{w} defines the connectivity from local point of view. That is, for any subset $S=\{u, v\}\subseteq V(G)$, let $\kappa_{G}(S)$ denote the maximum number of internally disjoint paths between $u$ and $v$ in $G$. Then $\kappa(G)=min\{\kappa_{G}(S)|S\subseteq V(G)$ and $|S|=2\}.$ As a generalization of the traditional connectivity, Chartrand {\em et al.}~\cite{c} introduced the {\it generalized $k$-connectivity} in $1984$. This parameter can measure the reliability of a network $G$ to connect any $k$ vertices in $G$. Let $S\subseteq V(G)$ and $\kappa_{G}(S)$ denote the maximum number $r$ of edge-disjoint trees $T_{1}, T_{2}, \ldots, T_{r}$ in $G$ such that $V(T_{i})\bigcap V(T_{j})=S$ for any $i, j \in \{1, 2, \ldots, r\}$ and $i\neq j$. For an integer $k$ with $2\leq k\leq n$, the {\it generalized $k$-connectivity} of a graph $G$ is defined as $\kappa_{k}(G)= min\{\kappa_{G}(S)|S\subseteq V(G)$ and $|S|=k\}$. The generalized $2$-connectivity is exactly the traditional connectivity. Li {\em et al.}~\cite{l4} derived that it is NP-complete for a general graph $G$ to decide whether there are $k$ internally disjoint trees connecting $S$, where $k$ is a fixed integer and $S\subseteq V(G).$ Some results~\cite{l2,l5} about the upper and lower bounds of the generalized connectivity are obtained. In addition, there are some results of generalized $3$-connectivity for some special kinds of graphs. For example, Chartrand {\em et al.}~\cite{ch} studied the generalized connectivity of complete graphs; Li {\em et al.}~\cite{LIS} characterized the minimally $2$-connected graphs $G$ with generalized connectivity $\kappa_{3}(G)=2$; Li {\em et al.}~\cite{l1} studied the generalized $3$-connectivity of Cartesian product graphs; Li {\em et al.}~\cite{l8} studied the generalized $3$-connectivity of graph products; Li {\em et al.}~\cite{l3} studied the generalized connectivity of the complete bipartite graphs; Li {\em et al.}~\cite{l6} studied the generalized $3$-connectivity of the star graphs and bubble-sort graphs and Li {\em et al.}~\cite{l7} studied the generalized $3$-connectivity of the Cayley graph generated by trees and cycles. For more results about the generalized connectivity, one can refer~\cite{lx, lh} for the detail.

In this paper, we study the generalized $3$-connectivity of some general $m$-regular and $m$-connected graphs $G_{n}$ constructed recursively and obtain that $\kappa_{3}(G_{n})=m-1$. As applications of the main result, the generalized $3$-connectivity of many famous networks such as the alternating group graph $AG_{n}$, the $k$-ary $n$-cube $Q_{n}^{k}$, the split-star network $S_{n}^{2}$ and the bubble-sort-star graph $BS_{n}$ etc. can be obtained directly.

The paper is organized as follows. In section 2, some notations and definitions are given. In section 3, the generalized $3$-connectivity of the regular graph $G_{n}$ is determined, which is the main result. In section $4$, as applications of the main result, the generalized $3$-connectivity of some famous networks are obtained. In section 5, the paper is concluded.
\section{Preliminary}

In this section, some terminologies and notation are introduced. For terminologies and notations undefined here, one can follow~\cite{B} for the detail.

\subsection{Terminologies and notation }
Let $G=(V, E)$ be a simple, undirected graph. Let $|V(G)|$ be the size of vertex set and $|E(G)|$ be the size of edge set. For a vertex $v$ in $G$, we denote by $N_{G}(v)$ the {\em neighbourhood} of the vertex $v$ in $G$ and $N_{G}[v]=N_{G}(v)\bigcup\{v\}$. Let $U \subseteq V(G)$, denote $N_G(U)=\bigcup\limits_{v\in U}N_{G}(v)-U$. Let $d_{G}(v)$ denote the degree of the vertex $v$ in $G$ and $\delta(G)$ denote the {\em minimum degree} of the graph $G$. The subgraph induced by $V^{\prime}$ in $G$, denoted by $G[V^{\prime}]$, is a graph whose vertex set is $V^{\prime}$ and the edge set is the set of all the edges of $G$ with both ends in $V^{\prime}$. A graph is said to be {\em $k$-regular} if for any vertex $v$ of $G$, $d_{G}(v)=k$. Two $xy$- paths $P$ and $Q$ in $G$ are {\em internally disjoint} if they have no common internal vertices, that is $V(P)\bigcap V(Q)=\{x, y\}$. Let $Y\subseteq V(G)$ and $X\subset V(G)\setminus Y$, the $(X, Y)$-paths is a family of internally disjoint paths starting at a vertex $x\in X$, ending at a vertex $y\in Y$ and whose internal vertices belong to neither $X$ nor $Y$. If $X=\{x\}$, the $(X, Y)$-paths is a family of internal disjoint paths whose starting vertex is $x$ and the terminal vertices are distinct in $Y$, which is referred to as a {\em $k$-fan} from $x$ to $Y$.

Let $[n]=\{1,2,\cdots,n\}$. Let $\Gamma$ be a finite group and $S$ be a subset of $\Gamma$, where the identity of the group does not belong to $S$. The {\em Cayley graph $Cay(\Gamma, S)$} is a digraph with vertex set $\Gamma$ and arc set $\{(g, g.s)| g\in \Gamma, s\in S\}$. If $S= S^{-1}$, then $Cay(\Gamma, S)$ is an undirected graph, where $S^{-1}=\{s^{-1}|s \in S\}$.

\begin{defi}\label{defi1} Let $n, r, a \geq 1$ be integers. Let $G_{n}$ be an  $n$-th regular graph, which can be recursively constructed as follows:
\begin{enumerate}
\item [{\rm (1)}] The $1$-th regular graph, say $G_{1}$, is a $r$-regular and $r$-connected graph with order $a$.

\item [{\rm (2)}] For $n\geq 2$, the $n$-th regular graph, say $G_{n}$, is a regular graph that consists of $p_{n}$ copies of $G_{n-1}$, say $G_{n-1}^{1}, G_{n-1}^{2}, \cdots, G_{n-1}^{p_{n}}$.

\item [{\rm (3)}] For each $u\in V(G_{n-1}^{i})$, it has two neighbours outside $G_{n-1}^{i}$, which are called the outside neighbours of $u$. In addition, the two outside neighbours of $u$ belong to different $G_{n-1}^{j}$s for $j\neq i$.

\item [{\rm (4)}] There are same edges between $G_{n-1}^{i}$ and $G_{n-1}^{j}$ for $i\neq j$. It can be checked that there are $\frac{2ap_{2}p_{3}\cdots p_{n-1}}{p_{n}-1}$ edges between $G_{n-1}^{i}$ and $G_{n-1}^{j}$ for $i\neq j$ and $i,j\in [p_{n}]$.

\item [{\rm (5)}] $\frac{2ap_{2}p_{3}\cdots p_{n-1}}{p_{n}-1}\geq r+2(n-2)+2$, where $r+2(n-2)\geq 4$.

\item [{\rm (6)}] $G_{n}$ is $m$-regular and $m$-connected, where $m=r+2(n-1)$.
\end{enumerate}
\end{defi}

For convenience, let $G_{n}=G_{n-1}^{1}\bigoplus G_{n-1}^{2}\bigoplus\cdots \bigoplus G_{n-1}^{p_{n}}$. By the definition of $G_{n}$, $|G_{n}|=N=ap_{2}p_{3}\cdots p_{n}$.

\subsection{Some networks can be regarded as the regular graph $G_{n}$}


\subsubsection{The alternating group graph $AG_{n}$}

The alternating group graph was introduced by Jwo $et$ $al.$~\cite{J} in $1993$. It is defined as follows.

\begin{defi}\label{defi6} Let $A_{n}$ be the alternating group of order $n$ with $n\geq 3$
and let $S=\{(12i), (1i2)|\\3\leq i\leq n\}$. The alternating group graph, denoted by $AG_{n}$, is defined as the Cayley graph $Cay(A_{n}, S)$.
\end{defi}

By the definition of $AG_{n}$, it is a $2(n-2)$-regular graph with $n!/2$ vertices. Let $A_{n}^{i}$ be the subset of $A_{n}$ that consists of all even permutations with element $i$ in the rightmost position and let $AG_{n-1}^{i}$ be the subgraph of $AG_{n}$ induced by $A_{n}^{i}$ for $i\in[n]$. Then $AG_{n-1}^{i}$ is isomorphic to $AG_{n-1}$ for each $i\in [n]$ and we call such an $AG_{n-1}^{i}$ a copy of $AG_{n-1}$. Thus, $AG_{n}$ can be decomposed into $n$ copies of $AG_{n-1}$, namely, $AG_{n-1}^{1}, AG_{n-1}^{2}, \cdots, AG_{n-1}^{n}$. For convenience, we denote $AG_{n}=AG_{n-1}^{1}\bigoplus AG_{n-1}^{2}\bigoplus\cdots \bigoplus AG_{n-1}^{n}$, where $\bigoplus$ just denotes the corresponding decomposition of $AG_{n}$. For each vertex $u\in V(AG_{n-1}^{i})$, it has $2(n-3)$ neighbours in $AG_{n-1}^{i}$ and two neighbours outside $AG_{n-1}^{i}$, which are called the outside neighbours of $u$.

\begin{lem}{\rm (\cite{ZHAO})}\label{lem4.9}
Let $AG_{n}=AG_{n-1}^{1}\bigoplus AG_{n-1}^{2}\bigoplus \ldots\bigoplus AG_{n-1}^{n}$ for $n\geq 3$. Then the following results hold.
\begin{enumerate}
\item [{\rm (1)}] For any vertex $u$ of $AG_{n-1}^{i}$, it has two outside neighbours.

\item [{\rm (2)}] For any copy $AG_{n-1}^{i}$, no two vertices in $AG_{n-1}^{i}$ have a common outside neighbour. In addition, $|N(AG_{n-1}^{i})|=(n-1)!$ and $|N(AG_{n-1}^{i})\bigcap V(AG_{n-1}^{j})|=(n-2)!$ for $i\neq j$.
\end{enumerate}
\end{lem}

\begin{lem}{\rm (\cite{J2, W1})}\label{lem4.10}
$\kappa(AG_{n})=2(n-2)$ for $n\geq 3$.
\end{lem}


\subsubsection{The $k$-ary $n$-cube network $Q_{n}^{k}$}

The $k$-ary $n$-cube network, denoted by $Q_{n}^{k}$, was introduced by S. Scott $et$ $al.$~\cite{sc} in $1994$. It is defined as follows.

\begin{defi}\label{defi7} The $k$-ary $n$-cube, denoted by $Q_{n}^{k}$, where $k\geq 2$ and $n\geq 1$ are integers, is a graph consisting of $k^{n}$ vertices, each of these vertices has the form $u=u_{n-1}u_{n-2}\cdots u_{0}$, where $u_{i}\in\{0,1,\cdots,k-1\}$ for $0\leq i\leq n-1$. Two vertices $u=u_{n-1}u_{n-2}\cdots u_{0}$ and $v=v_{n-1}v_{n-2}\cdots v_{0}$ in $Q_{n}^{k}$ are adjacent if and only if there exists an integer $j$, where $0\leq j\leq n-1$, such that $u_{j}=v_{j}\pm 1(mod $k$)$ and $u_{i}=v_{i}$ for every $i\in \{0,1,\cdots,k-1\}\setminus\{j\}$. In this case, $(u,v)$ is a $j$-dimensional edge.
\end{defi}

By the definition of $Q_{n}^{k}$, it is $2n$-regular for $k\geq 3$ and $n$-regular for $k=2$. Clearly, $Q_{1}^{k}$ is a cycle of length $k$ and $Q_{n}^{2}$ is the hypercube.

The $k$-ary $n$-cube $Q_{n}^{k}$ can be partitioned into $k$ disjoint subcubes along the $j$th-dimension for $j\in\{0, 1, 2, \cdots, n-1\}$, namely, $Q_{n-1}^{k}[0], Q_{n-1}^{k}[1], \cdots, Q_{n-1}^{k}[k-1]$. Then $Q_{n-1}^{k}[i]$ is isomorphic to the $k$-ary $(n-1)$-cube for $i\in\{0, 1, 2, \cdots, k-1\}$. For convenience, we denote $Q_{n}^{k}=Q_{n-1}^{k}[0]\bigoplus Q_{n-1}^{k}[1]\bigoplus\cdots \bigoplus Q_{n-1}^{k}[k-1]$, where $\bigoplus$ just denotes the corresponding decomposition of $Q_{n}^{k}$. For each vertex $u\in V(Q_{n-1}^{k}[i])$, it has $2n-2$ neighbours in $Q_{n-1}^{k}[i]$ and two neighbours outside $Q_{n-1}^{k}[i]$, which are called the outside neighbours of $u$.

\begin{lem}\label{lem4.11}
Let $Q_{n}^{k}=Q_{n-1}^{k}[0]\bigoplus Q_{n-1}^{k}[1]\bigoplus \ldots\bigoplus Q_{n-1}^{k}[k-1]$ for $k\geq 3$ and $n\geq 1$. Then the following results hold.
\begin{enumerate}
\item [{\rm (1)}] For any vertex $u$ of $Q_{n-1}^{k}[i]$, it has exactly two outside neighbours.
\item [{\rm (2)}] The outside neighbours of $u$ belong to different copies of $Q_{n-1}^{k}$. That is, no two vertices in $Q_{n-1}^{k}$ have a common outside neighbour.
\item [{\rm (3)}] $|N(Q_{n-1}^{k}[i])|=2k^{n-1}$ and $|N(Q_{n-1}^{k}[i])\bigcap V(Q_{n-1}^{k}[j])|=\frac{2k^{n-1}}{k-1}$ for $i\neq j$. That is, there are $\frac{2k^{n-1}}{k-1}$ independent crossed edges between two different $Q_{n-1}^{k}[i]$s.
\end{enumerate}
\end{lem}

\f {\bf Proof.} (1) Let $u=u_{1}u_{2}u_{3}\cdots u_{n-1}i\in V(Q_{n-1}^{k}[i])$. By definition $3$, $u^{\prime}=u_{1}u_{2}u_{3}\cdots u_{n-1}(i-1)$ and $u^{\prime\prime}=u_{1}u_{2}u_{3}\cdots u_{n-1}(i+1)$ are the two outside neighbours of $u$.

(2) Let $u=u_{1}u_{2}u_{3}\cdots u_{n-1}i\in V(Q_{n-1}^{k}[i])$. Then by (1), $u^{\prime}\in V(Q_{n-1}^{k}[i-1])$ and $u^{\prime\prime}\in V(Q_{n-1}^{k}[i+1])$ and they belong to different copies of $Q_{n-1}^{k}$.

(3) As any vertex in $Q_{n-1}^{k}[i]$ has two outside neighbours and $|Q_{n-1}^{k}[i]|=k^{n-1}$, then $|N(Q_{n-1}^{k}[i])|=2k^{n-1}$ and $|N(Q_{n-1}^{k}[i])\bigcap V(Q_{n-1}^{k}[j])|=\frac{2k^{n-1}}{k-1}$ for $i\neq j$.

\hfill\qed

\begin{lem}{\rm (\cite{D})}\label{lem4.12}
$\kappa(Q_{n}^{k})=2n$ for $k\geq 3$ and $n\geq 1$..
\end{lem}


\subsubsection{The split-star network $S_{n}^{2}$}

The split-star network, denoted by $S_{n}^{2}$, was proposed by E. Cheng $et$ $al.$~\cite{CH2} as an attractive variation of the star graph in $1998$. It is defined as follows.

\begin{defi}\label{defi8} Let $Sym(n)$ be symmetric group on $[n]$ and let $S=\{(1i)|2\leq i\leq n\}\bigcup\{(2i)|3\leq i\leq n\}$. The split-star network, denoted by $S_{n}^{2}$, is defined as the Cayley graph $Cay(Sym(n), S)$.
\end{defi}

By the definition of $S_{n}^{2}$, it is a $(2n-3)$-regular graph with $n!$ vertices. Let $V_{n}^{n:i}$ be the set of all vertices in $S_{n}^{2}$ with the $n$-th
position being $i$, that is, $V_{n}^{n:i}= \{u|u = u_{1}u_{2}\cdots u_{n-1}i\}$.
The set $\{V_{n}^{n:i}|1\leq i\leq n\}$ forms a partition of $V (S_{n}^{2})$. Let
$S_{n-1}^{2}[i]$ be the subgraph of $S_{n}^{2}$ induced by $V_{n}^{n:i}$. Then $S_{n-1}^{2}[i]$ is isomorphic to $S_{n-1}^{2}$ and we call such an $S_{n-1}^{2}[i]$ a copy of $S_{n-1}^{2}$. Thus, $S_{n}^{2}$ can be decomposed into $n$ copies of $S_{n-1}^{2}$, namely, $S_{n-1}^{2}[1], S_{n-1}^{2}[2], \cdots, S_{n-1}^{2}[n]$. For convenience, we denote $S_{n}^{2}=S_{n-1}^{2}[1]\bigoplus S_{n-1}^{2}[2]\bigoplus \ldots\bigoplus S_{n-1}^{2}[n]$, where $\bigoplus$ just denotes the corresponding decomposition of $S_{n}^{2}$.
For each vertex $u\in V(S_{n-1}^{2}[i])$, it has $2n-5$ neighbours in $S_{n-1}^{2}[i]$ and two neighbors outside $S_{n-1}^{2}[i]$, which are called outside neighbours of $u$.

\begin{lem}{\rm (\cite{CH})}\label{lem4.13}
Let $S_{n}^{2}=S_{n-1}^{2}[1]\bigoplus S_{n-1}^{2}[2]\bigoplus \ldots\bigoplus S_{n-1}^{2}[n]$ for $n\geq 3$. Then the following results hold.
\begin{enumerate}
\item [{\rm (1)}] For any vertex $u$ of $S_{n-1}^{2}[i]$ for $i\in[n]$, it has exactly two outside neighbours.
\item [{\rm (2)}] The outside neighbours of $u$ belong to different copies of $S_{n-1}^{2}$. That is, no two vertices in $S_{n-1}^{2}[i]$ have a common outside neighbour for $i\in[n]$.
\item [{\rm (3)}] $|N(S_{n-1}^{2}[i])|=2(n-1)!$ and $|N(S_{n-1}^{2}[i])\bigcap V(S_{n-1}^{2}[j])|=2(n-2)!$ for $i\neq j$. That is, there are $2(n-2)!$ independent crossed edges between two different $BS_{n-1}^{i}$s.

\end{enumerate}
\end{lem}

\begin{lem}{\rm (\cite{CH})}\label{lem4.14}
$\kappa(S_{n}^{2})=2n-3$ for $n\geq 3$.
\end{lem}


\subsubsection{The bubble-sort star graph $BS_{n}$}

The bubble-sort star graph, denoted by $BS_{n}$, was introduced by Z. Chou $et$ $al.$~\cite{chou} in $1996$. It is defined as follows.

\begin{defi}\label{defi9} Let $Sym(n)$ be symmetric group on $[n]$ and let $S=\{(1i)|2\leq i\leq n\}\bigcup\{(i,i+1)|2\leq i\leq n-1\}$. The $n$-dimensional bubble-sort star graph, denoted by $BS_{n}$, is defined as the Cayley graph $Cay(Sym(n), S)$.
\end{defi}

By the definition of $BS_{n}$, it is a $(2n-3)$-regular graph with $n!$ vertices. For an integer $i\in[n]$, let $BS_{n-1}^{i}$ be the graph induced by the vertex set $\{p_{1}p_{2}\cdots p_{n-1}i\}$, where $p_{1}p_{2}\cdots p_{n-1}$ ranges over all the permutations of $\{1,2,\cdots, i-1,i+1,\cdots,n\}$. Then $BS_{n-1}^{i}$ is isomorphic to $BS_{n-1}$ for each $i\in [n]$ and we call such an $BS_{n-1}^{i}$ a copy of $BS_{n-1}$. Thus, $BS_{n}$ can be decomposed into $n$ copies of $BS_{n-1}$, namely, $BS_{n-1}^{1}, BS_{n-1}^{2}, \cdots, BS_{n-1}^{n}$. For convenience, let $BS_{n}=BS_{n-1}^{1}\bigoplus BS_{n-1}^{2}\bigoplus \cdots BS_{n-1}^{n}$. For each vertex $u\in V(BS_{n-1}^{i})$, it has $2n-5$ neighbours in $BS_{n-1}^{i}$ and two neighbours outside $BS_{n-1}^{i}$, which are called the outside neighbours of $u$.

\begin{lem}{\rm(\cite{cai, wang})}\label{lem8} Let $BS_{n}=BS_{n-1}^{1}\bigoplus BS_{n-1}^{2}\bigoplus \ldots\bigoplus BS_{n-1}^{n}$, where $n\geq 4$. Then the following results hold.
\begin{enumerate}
\item [{\rm (1)}] For any vertex $u$ of $BS_{n-1}^{i}$ for $i\in[n]$, it has exactly two outside neighbours.

\item [{\rm (2)}] For any vertex $u$ of $BS_{n}$, the outside neighbours of $u$ belong to different copies of $BS_{n-1}$. That is, no two vertices in $BS_{n-1}^{i}$ have a common outside neighbour for $i\in \{1,2,\cdots, n\}$.
\item [{\rm (3)}] There are $2(n-2)!$ independent crossed edges between two different $BS_{n-1}^{i}$s.
\end{enumerate}
\end{lem}

\begin{lem}{\rm(\cite{cai})}\label{lem6}
$\kappa(BS_{n})=2n-3$ for $n\geq 3$.
\end{lem}


\section{The generalized $3$-connectivity of the regular graph $G_{n}$}

In this section, we will study the generalized $3$-connectivity of the regular graph $G_{n}$. The following lemmas are useful to our main result.

\begin{lem}{\rm(\cite{l5})}\label{lem1}
Let $G$ be a connected graph and $\delta$ be its minimum degree. Then $\kappa_{3}(G)\leq \delta$. Further, if there are two adjacent vertices of degree $\delta$, then $\kappa_{3}(G)\leq \delta-1$.
\end{lem}

\begin{lem}{\rm(\cite{l5})}\label{lem2}
Let $G$ be a connected graph with $n$ vertices. If $\kappa(G)=4k+r,$ where $k$ and $r$ are two integers with $k\geq 0$ and $r\in \{0, 1, 2, 3\},$ then $\kappa_{3}(G)\geq 3k+\lceil\frac{r}{2}\rceil$. Moreover, the lower bound is sharp.
\end{lem}

\begin{lem}{\rm(\cite{B})}\label{lem3}
Let $G=(V, E)$ be a $k$-connected graph, and let $X$ and $Y$ be subsets of $V(G)$ of cardinality at least $k$. Then there exists a family of $k$ pairwise disjoint $(X, Y)$-paths in $G$.
\end{lem}

\begin{lem}{\rm(\cite{B})}\label{lem4}
Let $G=(V, E)$ be a $k$-connected graph, let $x$ be a vertex of $G$, and let $Y\subseteq V\setminus \{x\}$ be a set of at least $k$ vertices of $G$. Then there exists a $k$-fan in $G$ from $x$ to $Y$, that is, there exists a family of $k$ internally disjoint $(x, Y)$-paths whose terminal vertices are distinct in $Y$.
\end{lem}

\begin{lem}\label{lem5}

Let $G_{n}$ and $r$ be the same as definition $1$. Let  $G_{n}=G_{n-1}^{1}\bigoplus G_{n-1}^{2} \bigoplus\ldots\bigoplus G_{n-1}^{p_{n}}$ and $H=G_{n-1}^{i_{1}}\bigoplus G_{n-1}^{i_{2}}\bigoplus \ldots\bigoplus G_{n-1}^{i_{l}}$ be the induced subgraph of $G_{n}$ on $\bigcup _{m=1}^{l}V(G_{n-1}^{i_{m}})$ for $2\leq l \leq p_{n}-1$, then $\kappa(H)\geq r+2(n-2)$, where $r+2(n-2)\geq 4$.
\end{lem}

\f {\bf Proof.} Without loss of generality, let $H=G_{n-1}^{1}\bigoplus G_{n-1 }^{2}\bigoplus \ldots\bigoplus G_{n-1}^{l}$. To prove the result, we just need to show that there are $r+2(n-2)$ internally disjoint paths for any two distinct vertices of $H$. Let $v_{1}, v_{2}\in V(H)$ and $v_{1}\neq v_{2}$, then the following two cases are considered.

Case 1. $v_{1}$ and $v_{2}$ belong to the same copy of $G_{n-1}$.

Without loss of generality, let $v_{1}, v_{2}\in V(G_{n-1}^{1})$. By definition $1(6)$, $\kappa(G_{n-1}^{1})=r+2(n-2)$. Then there are $r+2(n-2)$ internally disjoint paths between $v_{1}$ and $v_{2}$ in $G_{n-1}^{1}$.

Case 2. $v_{1}$ and $v_{2}$ belong to two different copies of $G_{n-1}$.

Without loss of generality, let $v_{1}\in V(G_{n-1}^{1})$ and $v_{2}\in V(G_{n-1}^{2})$. Select $r+2(n-2)$ vertices from $G_{n-1}^{1}\setminus\{v_{1}\}$, say $u_{1}, u_{2}, u_{3},\cdots,$ $u_{r+2(n-2)}$, such that the outside neighbour $u_{i}^{\prime}$ of $u_{i}$ belongs to $G_{n-1}^{2}\setminus\{v_{2}\}$ for each $i\in\{1,2,\cdots,r+2(n-2)\}$. By Definition $1$(5), this can be done. Let $S=\{u_{1}, u_{2}, u_{3},\cdots, u_{r+2(n-2)}\}$ and $S^{\prime}=\{u_{1}^{\prime}, u_{2}^{\prime}, u_{3}^{\prime},\cdots, u_{r+2(n-2)}^{\prime}\}$. By Definition $1(6)$, $\kappa(G_{n-1}^{1})=\kappa(G_{n-1}^{2})=r+2(n-2)$. By Lemma~\ref{lem4}, there exists a family of $r+2(n-2)$ internally disjoint $(v_{1}, S)$-paths $P_{1}, P_{2}, \cdots, P_{r+2(n-2)}$ such that the terminal vertex of $P_{i}$ is $u_{i}$. Similarly, there exists a family of $r+2(n-2)$ internally disjoint $(v_{2}, S^{\prime})$ paths $P_{1}^{\prime}, P_{2}^{\prime}, \cdots, P_{r+2(n-2)}^{\prime}$ such that the terminal vertex of $P_{i}^{\prime}$ is $u_{i}^{\prime}$. Let $\widehat{P_{i}}=P_{i}\bigcup u_{i}u_{i}^{\prime}\bigcup P_{i}^{\prime}$ for each $i\in\{1,2,\cdots, r+2(n-2)\}$, then $r+2(n-2)$ disjoint paths between $v_{1}$ and $v_{2}$ are obtained in $H$.\hfill\qed

\begin{lem}\label{lem6}
Let $G_{n}$ and $r$ be the same as definition $1$ and let $H=G_{n-1}^{i_{1}}\bigoplus G_{n-1}^{i_{2}}\bigoplus$
$G_{n-1}^{i_{3}}\bigoplus \cdots \bigoplus G_{n-1}^{i_{l}}$ be the induced subgraph of $G_{n}$ on $\bigcup _{j=1}^{l}V(G_{n-1}^{i_{j}})$ and $x\in V(H)$, where $l\geq 2$ and $n\geq 5$. If $d_{H}(x)=k$ and $Y\subseteq V(H)\setminus \{x\}$ with $|Y|=k$ such that $|Y \bigcap V(G_{n-1}^{i_{j}})|\leq r+2(n-2)$ for each $j\in\{1, 2, \cdots, l\}$. Then there exists a $k$-fan in $H$ from $x$ to $Y$.
\end{lem}

\f {\bf Proof.} Without loss of generality, let $H=G_{n-1}^{1}\bigoplus G_{n-1}^{2}\bigoplus G_{n-1}^{3}\bigoplus \cdots \bigoplus G_{n-1}^{l}$. Let $x\in V(H), d_{H}(x)=k$ and $Y\subseteq V(H)\setminus \{x\}$ with $|Y|=k$ such that $|Y \bigcap V(G_{n-1}^{j})|\leq r+2(n-2)$ for each $j\in\{1, 2, \cdots, l\}$. Clearly, $r+2(n-2)\leq k\leq r+2(n-1)$. To prove the result, the following three cases are considered.

Case $1$. $k=r+2(n-2)$.

By Lemma~\ref{lem5}, $\kappa(H)\geq r+2(n-2)$. By Lemma~\ref{lem4}, there exists a $[r+2(n-2)]$-fan in $H$ from $x$ to $Y$  and the result is desired.

Case $2$. $k=r+2(n-1)$.

Since $d_{H}(x)=r+2(n-1)$, then $V(H)$ contains the two outside neighbours $x^{\prime}$ and $x^{\prime\prime}$ of $x$. By definition $1(3)$, $x^{\prime}$ and $x^{\prime\prime}$ belong to different copies of $G_{n-1}$. Without loss of generality, let $x\in V(G_{n-1}^{1}), x^{\prime}\in V(G_{n-1}^{2})$ and $x^{\prime\prime}\in V(G_{n-1}^{3})$. Let $Y \bigcap V(G_{n-1}^{j})= A_{j}$ and $|A_{j}|=a_{j}$ for $1\leq j \leq l$. Then $a_{j}\leq r+2(n-2)$ and $\sum_{j=1}^{l}a_{j}=r+ 2(n-1)$. As $|Y|=r+2(n-1)$ and $|A_{1}|\leq r+2(n-2)$, there are at least two vertices of $Y$ outside $G_{n-1}^{1}$. We prove the result by considering $a_{j}$ for $j=2,3$ and the following two subcases are considered.

Subcase $2.1$. $a_{2}\geq 1$ and $a_{3}\geq 1$.

Let $a_{j}^{\prime}=a_{j}-1$ for $j=2, 3$ and $a_{j}^{\prime}=a_{j}$ for $j\in[l]\setminus\{2,3\}$. Then $\sum_{j=1}^{l}a_{j}^{\prime}=r+2(n-2)$. Now select $l-1$ pairwise disjoint vertex sets $M_{2}, M_{3}, \cdots, M_{l}$ in $G_{n-1}^{1}$ such that $|M_{j}|=a_{j}^{\prime}$ and for any vertex $v$ of $M_{j}$, one of the two outside neighbours of $v$ belongs to $G_{n-1}^{j}$ and $M_{j}\bigcap (A_{1}\bigcup \{x\})=\emptyset$ for $j\in \{ 2,3,\cdots, l\}$. By definition 1(5), this can be done. Let $M=A_{1}\bigcup M_{2}\bigcup \cdots \bigcup M_{l}$. As $|M|=r+2(n-2)$ and $\kappa(G_{n-1}^{1})=r+2(n-2)$. By Lemma~\ref{lem4}, there exist $l$ fans $F_{1}, F_{2}, \cdots, F_{l}$ in $G_{n-1}^{1}$ from $x$ to $A_{1}, M_{2}, \cdots, M_{l}$, respectively, where $F_{1}$ is a family of $a_{1}$ internally disjoint $(x, A_{1})$-paths whose terminal vertices are distinct in $A_{1}$ and $F_{j}$ is a family of $a_{j}^{\prime}$ internally disjoint $(x, M_{j})$-paths whose terminal vertices are distinct in $M_{j}$ for $2\leq j\leq l$. See Figure $1$.
\begin{figure}[!ht]
\begin{center}
\vskip0.1cm
\includegraphics[scale=0.8]{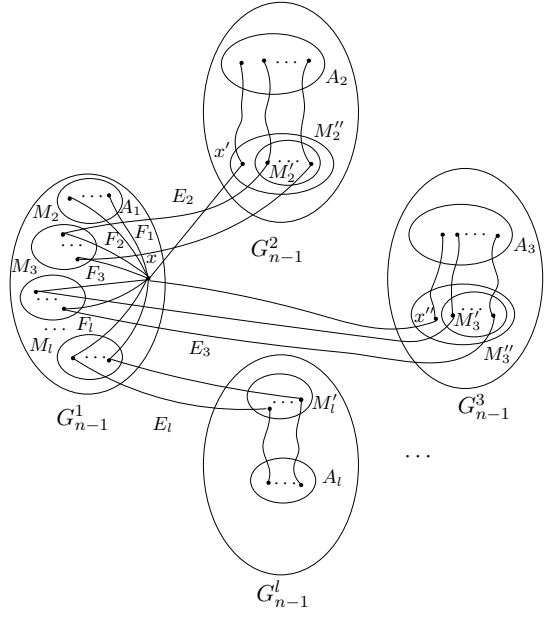}
\end{center}
\vskip0.1cm
\caption{Illustration of Subcase $2.1$ for $A_{j0}=\emptyset$ for each $j\in\{2,3\cdots,l\}$}
\end{figure}
Let $M_{j}^{\prime}=\{y^{\prime}|y^{\prime}$ is the outside neighbour of $y$ such that $y^{\prime}\in V(G_{n-1}^{j})$ for each $y\in M_{j}\}$ and $E_{j}=\{yy^{\prime}\in E(G_{n})| y\in M_{j}$ and $y^{\prime} \in M_{j}^{\prime}\}$ for $2\leq j\leq l$. Let $M_{2}^{\prime\prime}=M_{2}^{\prime}\bigcup \{x^{\prime}\}$ and $M_{3}^{\prime\prime}=M_{3}^{\prime}\bigcup \{x^{\prime\prime}\}$, then $|M_{2}^{\prime\prime}|=a_{2}$ and $|M_{3}^{\prime\prime}|=a_{3}$. Let $M_{j}^{\prime\prime}\bigcap A_{j}=A_{j0}$ for $j=2,3$ and $M_{j}^{\prime}\bigcap A_{j}=A_{j0}$ for $4\leq j\leq l$. Let $M_{j}^{\prime\prime}\setminus A_{j0}=A_{j1}$ for $j=2,3$ and $M_{j}^{\prime}\setminus A_{j0}=A_{j1}$ for $4\leq j\leq l$, and let $A_{j}\setminus A_{j0}=A_{j2}$ for $2\leq j\leq l$. Then $|A_{j1}|=|A_{j2}|=a_{j}-|A_{j0}|$ for $2\leq j\leq l$. By definition $1(6)$, $\kappa(G_{n-1}^{j})=r+2(n-2)$. As $\kappa(G_{n-1}^{j}\setminus A_{j0})\geq r+2(n-2)-|A_{j0}|\geq a_{j}-|A_{j0}|$. By Lemma~\ref{lem3}, there exists a family of $a_{j}-|A_{j0}|$ pairwise disjoint $(A_{j1}, A_{j2})$-paths $F_{j}^{\prime}$ in $G_{n-1}^{j}$ for $2\leq j\leq l$.

Finally, by combining the $l$ fans $F_{1}, F_{2}, \cdots, F_{l}$, the edge sets $E_{2}, \cdots, E_{l}$, the edges $xx^{\prime}, xx^{\prime\prime}$ and the paths $F_{2}^{\prime}, \cdots, F_{l}^{\prime}$, we can obtain a $[r+2(n-1)]$-fan from $x$ to $Y$ in $H$.

Subcase $2.2$. At least one of $a_{2}, a_{3}=0$.

Without loss of generality, we assume $a_{2}=0$ and the following three subcases are considered.

Subcase $2.2.1$. $a_{2}=0$ and $a_{3}\geq 2$.

Since $a_{2}=0$ and $a_{3}\geq 2$, see Figure $2$. Let $a_{j}^{\prime}=a_{j}-2$ for $j=3$ and $a_{j}^{\prime}=a_{j}$ for $j\in[l]\setminus\{3\}$. Then select $l-2$ pairwise disjoint vertex sets $M_{3}, M_{4}, \cdots, M_{l}$ in $G_{n-1}^{1}$ such that $|M_{j}|=a_{j}^{\prime}$ and for any vertex $v$ of $M_{j}$, one of the two outside neighbours of $v$ belongs to $G_{n-1}^{j}$ and $M_{j}\bigcap (A_{1}\bigcup \{x\})=\emptyset$ for each $j\in \{ 3,4,\cdots, l\}$. By definition 1(5), this can be done. Let $M=A_{1}\bigcup M_{3}\bigcup \cdots \bigcup M_{l}$. As $|M|=r+2(n-2)$ and $\kappa(G_{n-1}^{1})=r+2(n-2)$ by definition $1(6)$. By Lemma~\ref{lem4}, there exist $l-1$ fans $F_{1}, F_{3}, \cdots, F_{l}$ in $G_{n-1}^{1}$ from $x$ to $A_{1}, M_{3}, \cdots, M_{l}$, respectively.
\begin{figure}[!ht]
\begin{center}
\vskip0.1cm
\includegraphics[scale=0.6]{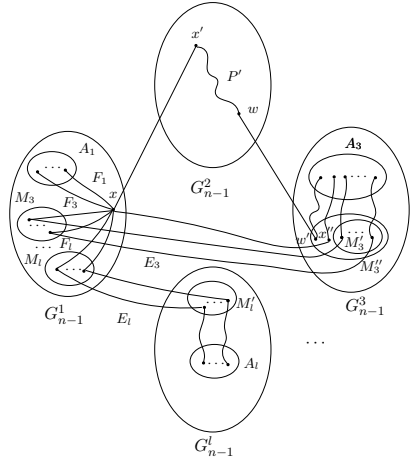}
\end{center}
\vskip0.1cm
\caption{Illustration of Subcase $2.2.1$ for $A_{j0}=\emptyset$ for each $j\in\{3,4\cdots,l\}$}
\end{figure}
Let $M_{j}^{\prime}=\{y^{\prime}|y^{\prime}$ is the outside neighbour of $y$ such that $y^{\prime}\in V(G_{n-1}^{j})$ for each $y\in M_{j}\}$ and $E_{j}=\{yy^{\prime}\in E(G_{n})| y\in M_{j}$ and $y^{\prime} \in M_{j}^{\prime}\}$ for $3\leq j\leq l$. Let $w\in V(G_{n-1}^{2})$ and one of the outside neighbours $w^{\prime}$ of $w$ belongs to $V(G_{n-1}^{3})$ and $w^{\prime}\notin \{x^{\prime\prime}\} \bigcup M_{3}^{\prime}$. By definition $1(5)$, this can be done. Then there exists a path $P^{\prime}$ between $x^{\prime}$ and $w$. Let $M_{3}^{\prime\prime}=M_{3}^{\prime}\bigcup \{x^{\prime\prime}, w^{\prime}\}$, then $|M_{3}^{\prime\prime}|=a_{3}$. Let $M_{j}^{\prime\prime}\bigcap A_{j}=A_{j0}$ for $j=3$ and $M_{j}^{\prime}\bigcap A_{j}=A_{j0}$ for $4\leq j\leq l$. Let $M_{j}^{\prime\prime}\setminus A_{j0}=A_{j1}$ for $j=3$ and $M_{j}^{\prime}\setminus A_{j0}=A_{j1}$ for $4\leq j\leq l$, and let $A_{j}\setminus A_{j0}=A_{j2}$ for $3\leq j\leq l$. Then $|A_{j1}|=|A_{j2}|=a_{j}-|A_{j0}|$ for $3\leq j\leq l$. By definition $1(6)$, $\kappa(G_{n-1}^{j})=r+2(n-2)$. As $\kappa(G_{n-1}^{j}\setminus A_{j0})\geq r+2(n-2)-|A_{j0}|\geq a_{j}-|A_{j0}|$. By Lemma~\ref{lem3}, there exists a family of $a_{j}-|A_{j0}|$ pairwise disjoint $(A_{j1}, A_{j2})$-paths $F_{j}^{\prime}$ in $AG_{n-1}^{j}$ for $3\leq j\leq l$.

Next, by combining the $l-1$ fans $F_{1}, F_{3} \cdots, F_{l}$, the edge sets $E_{3}, E_{4}, \cdots, E_{l}$, the edges $xx^{\prime}, xx^{\prime\prime}, ww^{\prime}$, the path $P^{\prime}$ and the paths $F_{3}^{\prime}, \cdots, F_{l}^{\prime}$, we can obtain a $[r+2(n-1)]$-fan from $x$ to $Y$ in $H$.

Subcase $2.2.2$. $a_{2}=0$ and $a_{3}=1$.

Since $a_{2}=0$ and $a_{3}=1$, there must exist a part $G_{n-1}^{k}$ such that $a_{k}\geq 1$ for $k\in\{4,5,\cdots,l\}$. Let $a_{j}^{\prime}=a_{j}-1$ for $j=3, k$ and $a_{j}^{\prime}=a_{j}$ for $j\in[l]\setminus\{3, k\}$.

Then select $l-2$ pairwise disjoint vertex sets $M_{3}, M_{4}, \cdots, M_{l}$ in $G_{n-1}^{1}$ such that $|M_{j}|=a_{j}^{\prime}$ and for any vertex $v$ of $M_{j}$, one of the two outside neighbours of $v$ belongs to $G_{n-1}^{j}$ and $M_{j}\bigcap (A_{1}\bigcup \{x\})=\emptyset$ for each $j\in \{ 3,4,\cdots, l\}$. Let $M=A_{1}\bigcup M_{3}\bigcup \cdots \bigcup M_{l}$. As $|M|=r+2(n-2)$ and $\kappa(G_{n-1}^{1})=r+2(n-3)$ by definition $1(6)$. By Lemma~\ref{lem4}, there exist $l-1$ fans $F_{1}, F_{3}, \cdots, F_{l}$ in $G_{n-1}^{1}$ from $x$ to $M$, where $F_{j}$ is a family of $a_{j}^{\prime}$ internally disjoint $(x, M_{j})$-paths whose terminal vertices are distinct in $M_{j}$ for $3\leq j\leq l$.

Let $M_{j}^{\prime}=\{y^{\prime}|y^{\prime}$ is the outside neighbour of $y$ such that $y^{\prime}\in V(G_{n-1}^{j})$ for each $y\in M_{j}\}$ and $E_{j}=\{yy^{\prime}\in E(G_{n})| y\in M_{j}$ and $y^{\prime} \in M_{j}^{\prime}\}$ for $3\leq j\leq l$. Let $w\in V(G_{n-1}^{2})$ such that one of the outside neighbours $w^{\prime}$ of $w$ belongs to $G_{n-1}^{k}$ and $w^{\prime}\notin M_{k}^{\prime}$. Then there exists a path $P^{\prime}$ from $x^{\prime}$ to $w$ in $G_{n-1}^{2}$. Let $M_{k}^{\prime\prime}=M_{k}^{\prime}\bigcup \{w^{\prime}\}$ and $M_{3}^{\prime\prime}=M_{3}^{\prime}\bigcup \{x^{\prime\prime}\}$, then $|M_{k}^{\prime\prime}|=a_{k}$ and $|M_{3}^{\prime\prime}|=a_{3}$.
Then prove the result similar as Subcase $2.1$, we can obtain a $[r+2(n-1)]$-fan from $x$ to $Y$ in $H$.

Subcase $2.2.3$. $a_{2}=0$ and $a_{3}=0$.

In this case, there exists a part $G_{n-1}^{k}$ such that $a_{k}\geq 2$ for $k\in \{4, 5, \cdots, l\}$ or there exist two parts $G_{n-1}^{i}$ and $G_{n-1}^{m}$ such that $a_{i}, a_{m}\geq 1$ for $i, m\in \{4, 5, \cdots, l\}$.

Subcase $2.2.3.1$. There exists a part $G_{n-1}^{k}$ such that $a_{k}\geq 2$ for $k\in \{4, 5, \cdots, l\}$.

For this case, see Figure 3. Let $a_{j}^{\prime}=a_{j}-2$ for $j=k$ and $a_{j}^{\prime}=a_{j}$ for $j\neq k$. Then select $l-3$ pairwise disjoint vertex sets $M_{4}, M_{5}, \cdots, M_{l}$ in $G_{n-1}^{1}$ such that $|M_{j}|=a_{j}^{\prime}$ and for any vertex $v$ of $M_{j}$, one of the two outside neighbours of $v$ belongs to $G_{n-1}^{j}$ and $M_{j}\bigcap (A_{1}\bigcup \{x\})=\emptyset$ for each $j\in \{4,\cdots, l\}$. Let $M=A_{1}\bigcup M_{4}\bigcup \cdots \bigcup M_{l}$. As $|M|=r+2(n-2)$ and $\kappa(G_{n-1}^{1})=r+2(n-2)$ by definition $1(6)$. By Lemma~\ref{lem4}, there exist $l-2$ fans $F_{1}, F_{4}, \cdots, F_{l}$ in $G_{n-1}^{1}$ from $x$ to $M$, where $F_{j}$ is a family of $a_{j}^{\prime}$ internally disjoint $(x, M_{j})$-paths whose terminal vertices are distinct in $M_{j}$ for $4\leq j\leq l$. Let $M_{j}^{\prime}=\{y^{\prime}|y^{\prime}$ is the outside neighbour of $y$ such that $y^{\prime}\in V(G_{n-1}^{j})$ for each $y\in M_{j}\}$ and $E_{j}=\{yy^{\prime}\in E(G_{n})| y\in M_{j}$ and $y^{\prime} \in M_{j}^{\prime}\}$ for $4\leq j\leq l$. Let $u\in V(G_{n-1}^{2})$ and one of the outside neighbours $u^{\prime}$ of $u$ belongs to $V(G_{n-1}^{k})$ and $u^{\prime}\notin M_{k}^{\prime}$. Let $v\in V(G_{n-1}^{3})$ and one of the outside neighbours $v^{\prime}$ of $v$ belongs to $V(G_{n-1}^{k})$ and $v^{\prime}\notin \{u^{\prime}\}\bigcup M_{k}^{\prime}$. Then there exists a path $P_{1}$ between $x^{\prime}$ and $u$ in $G_{n-1}^{2}$ and a path $P_{2}$ between $x^{\prime\prime}$ and $v$ in $G_{n-1}^{3}$. Let $M_{k}^{\prime\prime}=M_{k}^{\prime}\bigcup \{u^{\prime}, v^{\prime}\}$, then $|M_{k}^{\prime\prime}|=a_{k}$. Then prove the result similar as Subcase $2.2.1$, we can obtain a $[r+2(n-1)]$-fan from $x$ to $Y$ in $H$.

\begin{figure}[!ht]
\begin{center}
\vskip0.1cm
\includegraphics[scale=0.6]{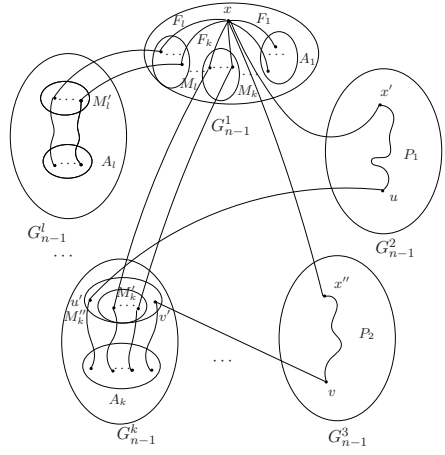}
\end{center}
\vskip0.1cm
\caption{Illustration of Subcase $2.2.3.1$}
\end{figure}
Subcase $2.2.3.2$. There exists two parts $G_{n-1}^{i}$ and $G_{n-1}^{m}$ such that $a_{i}, a_{m}\geq 1$ for $i, m\in \{4, 5, \cdots, l\}$.

For this case, see Figure 4. Let $a_{j}^{\prime}=a_{j}-1$ for $j=i,m$ and $a_{j}^{\prime}=a_{j}$ for $j\neq i,m$. Then select $l-3$ pairwise disjoint vertex sets $M_{4}, M_{5}, \cdots, M_{l}$ in $G_{n-1}^{1}$ such that $|M_{j}|=a_{j}^{\prime}$ and for any vertex $v$ of $M_{j}$, one of the two outside neighbours of $v$ belongs to $G_{n-1}^{j}$ and $M_{j}\bigcap (A_{1}\bigcup \{x\})=\emptyset$ for each $j\in \{4,\cdots, l\}$. Let $M=A_{1}\bigcup M_{4}\bigcup \cdots \bigcup M_{l}$. As $|M|=r+2(n-2)$ and $\kappa(G_{n-1}^{1})=r+2(n-2)$ by definition $1(6)$. By Lemma~\ref{lem4}, there exist $l-2$ fans $F_{1}, F_{4}, \cdots, F_{l}$ in $G_{n-1}^{1}$ from $x$ to $M$, where $F_{j}$ is a family of $a_{j}^{\prime}$ internally disjoint $(x, M_{j})$-paths whose terminal vertices are distinct in $M_{j}$ for $4\leq j\leq l$. Let $M_{j}^{\prime}=\{y^{\prime}|y^{\prime}$ is the outside neighbour of $y$ such that $y^{\prime}\in V(G_{n-1}^{j})$ for each $y\in M_{j}\}$ and $E_{j}=\{yy^{\prime}\in E(G_{n})| y\in M_{j}$ and $y^{\prime} \in M_{j}^{\prime}\}$ for $4\leq j\leq l$. Let $u\in V(G_{n-1}^{2})$ and one of the outside neighbours $u^{\prime}$ of $u$ belongs to $V(G_{n-1}^{i})$ and $u^{\prime}\notin M_{i}^{\prime}$. Let $v\in V(G_{n-1}^{3})$ and one of the outside neighbours $v^{\prime}$ of $v$ belongs to $V(G_{n-1}^{m})$ and $v^{\prime}\notin M_{m}^{\prime}$. Then there exists a path $P_{1}$ between $x^{\prime}$ and $u$ in $G_{n-1}^{2}$ and a path $P_{2}$ between $x^{\prime\prime}$ and $v$ in $G_{n-1}^{3}$. Let $M_{i}^{\prime\prime}=M_{i}^{\prime}\bigcup \{u^{\prime}\}$ and $M_{m}^{\prime\prime}=M_{m}^{\prime}\bigcup \{v^{\prime}\}$, then $|M_{i}^{\prime\prime}|=a_{i}$ and $|M_{m}^{\prime\prime}|=a_{m}$. Then prove the result similar as Subcase $2.2.1$, we can obtain a $[r+2(n-1)]$-fan from $x$ to $Y$ in $H$.
\begin{figure}[!ht]
\begin{center}
\vskip0.1cm
\includegraphics[scale=0.6]{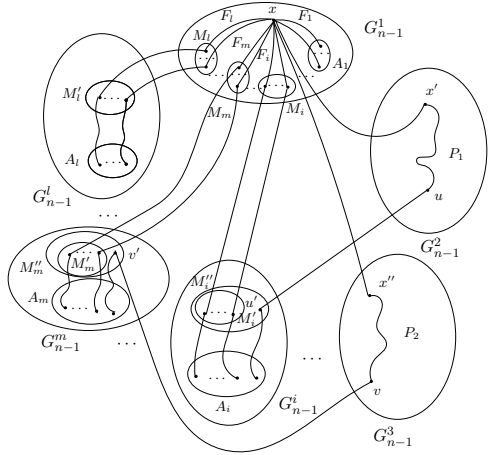}
\end{center}
\vskip0.1cm
\caption{Illustration of Subcase $2.2.3.2$}
\end{figure}

Case $3$. $k=r+2n-3$.

Since $d_{H}(x)=r+2n-3$, then $V(H)$ contains one outside neighbour of $x$. Then prove the result similar as Case $2$. To avoid duplication of discussions, the discussions for this case is omitted.

Hence, there exists a $k$-fan from $x$ to $Y$ in $H$ in any case.
\hfill\qed

\begin{lem}\label{lem7}
Let $G_{n}$ and $r$ be the same as definition $1$ and $H=G_{n-1}^{i}\bigoplus G_{n-1}^{j}$ for $i\neq j$. If $x\in V(G_{n-1}^{i}), y\in V(G_{n-1}^{j})$ and $d_{H}(x)=d_{H}(y)=r+2n-3$, then there exist $r+2n-3$ internally disjoint paths between $x$ and $y$ in $H$.
\end{lem}

\f {\bf Proof.} Without loss of generality, let $H=G_{n-1}^{1}\bigoplus G_{n-1}^{2}$, $x\in V(G_{n-1}^{1}), y\in V(G_{n-1}^{2})$ and $d_{H}(x)=d_{H}(y)=r+2n-3$. To prove the main result, the following two cases are considered.

Case 1. $x$ and $y$ are not adjacent.

Let $Y=N_{H}(y)=\{y_{1}, y_{2}, \cdots, y_{r+2n-3}\}$, then $x \notin Y$. Otherwise, $x$ and $y$ are adjacent. Clearly, $|Y \bigcap V(G_{n-1}^{m})|\leq r+2(n-2)$ for $m=1, 2$ and $|Y|=r+2n-3$. By Lemma~\ref{lem6}, there exist $r+2n-3$ internally disjoint paths $P_{1}, P_{2}, \cdots, P_{r+2n-3}$ in $H$ from $x$ to $Y$ whose terminal vertices are distinct in $Y$. If none of the paths $P_{i}$s for $1\leq i\leq r+2n-3$ contains $y$ as an internal vertex, then combining the edges from $y$ to $Y$ and the paths $P_{i}$s for $1\leq i\leq r+2n-3$, $r+2n-3$ internally disjoint paths between $x$ and $y$ in $H$ can be obtained. If not, there exists only one path which contains $y$ as an internal vertex as $P_{i}$s for $1\leq i\leq r+2n-3$ are internally disjoint. Assume that $P_{1}$ contains $y$ as an internal vertex and the terminal vertex of $P_{1}$ is $y_{1}$. Then $P_{1}$ contains a subpath $\widetilde{P}_{1}$ from $x$ to $y$. Then combining the edges from $y$ to $Y\setminus \{y_{1}\}, \widetilde{P}_{1}$ and the paths $P_{i}$s for $2\leq i\leq r+2n-3$, $r+2n-3$ internally disjoint $(x,y)$-paths in $H$ can be obtained.

Case 2. $x$ and $y$ are adjacent.

Choose $r+2(n-2)$ vertices $x_{1}, x_{2}, \cdots, x_{r+2(n-2)}$ from $G_{n-1}^{1}\setminus\{x\}$ such that one of the outside neighbours of $x_{i}$ belongs to $G_{n-1}^{2}\setminus\{y\}$ for each $i\in \{1,2,\cdots, r+2(n-2)\}$. Let $X=\{x_{1}, x_{2}, \cdots, x_{r+2(n-2)}\}$ and $X^{\prime}=\{x_{i}^{\prime}|x_{i}^{\prime}$ is the outside neighbour of $x_{i}$ and $x_{i}^{\prime}\in V(G_{n-1}^{2})\}$. By definition 1(5), this can be done. By definition $1(6)$, $\kappa(G_{n-1}^{1})=\kappa (G_{n-1}^{2})=r+2(n-2)$. By Lemma~\ref{lem4}, there exist $r+2(n-2)$ internally disjoint paths $P_{1}, P_{2}, \cdots, P_{r+2(n-2)}$ from $x$ to $X$ such that the terminal vertex of $P_{i}$ is $x_{i}$ in $G_{n-1}^{1}$ and $r+2(n-2)$ internally disjoint paths $P_{1}^{\prime}, P_{2}^{\prime}, \cdots, P_{r+2(n-2)}^{\prime}$ from $y$ to $X^{\prime}$ such that the terminal vertex of $P_{i}^{\prime}$ is $x_{i}^{\prime}$ in $G_{n-1}^{2}$ for each $i\in\{1,2,\cdots, r+2(n-2)\}$. Let $\widetilde{P}_{r+2n-3}=xy$ and $\widetilde{P_{i}}=xP_{i}x_{i}x_{i}^{\prime}P_{i}^{\prime}y$ for $1\leq i\leq r+2(n-2)$. Then $r+2n-3$ internally disjoint paths $\widetilde{P_{i}}$s for $1\leq i\leq r+2n-3$ between $x$ and $y$ in $H$ are obtained.
\hfill\qed

\begin{lem}\label{lem8}
Let $G_{n}$ and $r$ be the same as definition $1$, $G_{n}=G_{n-1}^{1}\bigoplus G_{n-1}^{2} \bigoplus \ldots \bigoplus G_{n-1}^{p_{n}}$ and $S=\{v_{1}, v_{2}, v_{3}\}$, where $v_{1}, v_{2}$ and $v_{3}$ are any three distinct vertices of $V(G_{n-1}^{i})$ for some $i$ with $1\leq i\leq p_{n}$. If there exist $r+2n-5$ internally disjoint trees connecting $S$ in $G_{n-1}^{i}$, then there exist $r+2n-3$ internally disjoint trees connecting $S$ in $G_{n}$.
\end{lem}

\f {\bf Proof.} Without loss of generality, let $S\subseteq V(G_{n-1}^{1})$. Note that there exist $r+2n-5$ internally disjoint trees $T_{1}, T_{2}, \ldots, T_{r+2n-5}$ connecting $S$ in $G_{n-1}^{1}$. As $v_{i}$ has two outside neighbours $v_{i}^{\prime}$ and $v_{i}^{\prime\prime}$ for each $i\in \{1, 2, 3\}$ and any two distinct vertices of $G_{n-1}^{1}$ have different outside neighbours by definition $1(3)$. Hence, $M=\{v_{1}^{\prime}, v_{2}^{\prime}, v_{3}^{\prime}, v_{1}^{\prime\prime}, v_{2}^{\prime\prime}, v_{3}^{\prime\prime}\}$ contains exactly $6$ distinct vertices. In addition, each copy of $G_{n-1}$ contains at most three vertices of them. To prove the result, the following three cases are considered.

Case $1$. There exists a copy of $G_{n-1}$ which contains three vertices of $M$.

Without loss of generality, let $\{v_{1}^{\prime}, v_{2}^{\prime}, v_{3}^{\prime}\}\subseteq V(G_{n-1}^{2})$ and $\{v_{1}^{\prime\prime}, v_{2}^{\prime\prime}, v_{3}^{\prime\prime}\}\subseteq \bigcup_{i=3}^{p_{n}}V(G_{n-1}^{i})$. As $G_{n-1}^{2}$ and $G_{n}[\bigcup _{i=3}^{p_{n}}V(G_{n-1}^{i})]$ as subgraphs of $G_{n}$ are both connected, there is a tree, say $T_{r+2n-4}^{\prime}$, connecting $v_{1}^{\prime}$, $v_{2}^{\prime}$ and $v_{3}^{\prime}$ in $G_{n-1}^{2}$ and a tree, say $T_{r+2n-3}^{\prime}$, connecting $v_{1}^{\prime\prime}, v_{2}^{\prime\prime}$ and $v_{3}^{\prime\prime}$ in $G_{n}[\bigcup _{i=3}^{p_{n}}V(G_{n-1}^{i})]$, respectively. Let $T_{r+2n-4}=T_{r+2n-4}^{\prime}\bigcup v_{1}v_{1}^{\prime}\bigcup v_{2}v_{2}^{\prime}\bigcup v_{3}v_{3}^{\prime}$ and $T_{r+2n-3}=T_{r+2n-3}^{\prime}\bigcup v_{1}v_{1}^{\prime\prime}\bigcup v_{2}v_{2}^{\prime\prime}\bigcup v_{3}v_{3}^{\prime\prime}$, then $V(T_{r+2n-4})\bigcap V(T_{r+2n-3})=S$. Combining the trees $T_{i}$s for $1\leq i\leq r+2n-3$, $r+2n-3$ internally disjoint trees connecting $S$ are obtained in $G_{n}$.

Case $2$. There exists a copy of $G_{n-1}$ which contains two vertices of $M$ and all other copies of $G_{n-1}$ contain at most two vertices of $M$.

Without loss of generality, let $v_{1}^{\prime}, v_{2}^{\prime}\in V(G_{n-1}^{2})$ and $v_{3}^{\prime}\in V(G_{n-1}^{3})$. The following two subcases are considered.

Subcase $2.1$. $G_{n-1}^{3}$ contains only the vertex $v_{3}^{\prime}$ of $M\setminus \{v_{1}^{\prime}, v_{2}^{\prime}\}$.

As $G_{n}[\bigcup_{i=2}^{3}V(G_{n-1}^{i})]$ and $G_{n}[\bigcup_{i=4}^{p_{n}}V(G_{n-1}^{i})]$ as subgraphs of $G_{n}$ are both connected, there is a tree, say $T_{r+2n-4}^{\prime}$, connecting $v_{1}^{\prime}, v_{2}^{\prime}$ and $v_{3}^{\prime}$ in $G_{n}[\bigcup_{i=2}^{3}V(G_{n-1}^{i})]$ and a tree, say $T_{r+2n-3}^{\prime}$, connecting $v_{1}^{\prime\prime}, v_{2}^{\prime\prime}$ and $v_{3}^{\prime\prime}$ in $G_{n}[\bigcup_{i=4}^{p_{n}}V(G_{n-1}^{i})]$, respectively. Let $T_{r+2n-4}=T_{r+2n-4}^{\prime}\bigcup v_{1}v_{1}^{\prime}\bigcup v_{2}v_{2}^{\prime} \bigcup v_{3}v_{3}^{\prime}$ and $T_{r+2n-3}=T_{r+2n-3}^{\prime}\bigcup v_{1}v_{1}^{\prime\prime}\bigcup v_{2}v_{2}^{\prime\prime}\bigcup v_{3}v_{3}^{\prime\prime}$, then $V(T_{r+2n-4})\\ \bigcap V(T_{r+2n-3})=S$. Combining the trees $T_{i}$s for $1\leq i\leq r+2n-3$, $r+2n-3$ internally disjoint trees connecting $S$ are obtained in $G_{n}$.

Subcase $2.2$. $G_{n-1}^{3}$ contains the vertex $v_{3}^{\prime}$ and a vertex of $M\setminus \{v_{1}^{\prime}, v_{2}^{\prime}, v_{3}^{\prime}\}$.

Without loss of generality, let $v_{3}^{\prime}, v_{1}^{\prime\prime}\in V(G_{n-1}^{3})$ and the following two subcases are considered.

Subcase $2.2.1$. $v_{3}^{\prime\prime}$ and $v_{2}^{\prime\prime}$ belong to different copies of $G_{n-1}$.

Without loss of generality, let $v_{3}^{\prime\prime}\in V(G_{n-1}^{4})$ and $v_{2}^{\prime\prime}\in V(G_{n-1}^{5})$. As $G_{n}[V(G_{n-1}^{2})\bigcup V(G_{n-1}^{4}\\)]$ is connected, there is a tree, say $T_{r+2n-4}^{\prime}$, connecting $v_{1}^{\prime}, v_{2}^{\prime}$ and $v_{3}^{\prime\prime}$ in $G_{n}[V(G_{n-1}^{2})\bigcup V(G_{n-1}^{4})]$. In addition, there is a tree, say $T_{r+2n-3}^{\prime}$, connecting $v_{1}^{\prime\prime}, v_{2}^{\prime\prime}$ and $v_{3}^{\prime}$ in $G_{n}[\bigcup_{i\in[p_{n}]\setminus\{1,2,4\}}V(G_{n-1}^{i})]$ as it is connected. Let $T_{r+2n-4}=T_{r+2n-4}^{\prime}\bigcup v_{1}v_{1}^{\prime}\bigcup v_{2}v_{2}^{\prime}\bigcup v_{3}v_{3}^{\prime\prime}$ and $T_{r+2n-3}=T_{r+2n-3}^{\prime}\bigcup v_{1}v_{1}^{\prime\prime}\\ \bigcup v_{2}v_{2}^{\prime\prime}\bigcup v_{3}v_{3}^{\prime}$, then $V(T_{r+2n-4})\bigcap V(T_{r+2n-3})=S$. Combining the trees $T_{i}$s for $1\leq i\leq r+2n-3$, $r+2n-3$ internally disjoint trees connecting $S$ are obtained in $G_{n}$.

Subcase $2.2.2$. $v_{3}^{\prime\prime}$ and $v_{2}^{\prime\prime}$ belong to the same copy of $G_{n-1}$.

Without loss of generality, let $v_{3}^{\prime\prime}, v_{2}^{\prime\prime}\in V(G_{n-1}^{4})$. As $v_{3}$ is one of the outside neighbours of $v_{3}^{\prime}$ and it has exactly two outside neighbours. Then let the other outside neighbour of $v_{3}^{\prime}$ be $u$. If $u\notin V(G_{n-1}^{4})$, then $G_{n}[\bigcup _{i\in[p_{n}]\setminus\{1,3,4\}}V(G_{n-1}^{i})]$ contains a tree $T_{r+2n-4}^{\prime}$ connecting $v_{1}^{\prime}, v_{2}^{\prime}$ and $u$. Let $T_{r+2n-4}=T_{r+2n-4}^{\prime}\bigcup v_{1}v_{1}^{\prime}\bigcup v_{2}v_{2}^{\prime}\bigcup v_{3}v_{3}^{\prime}\bigcup v_{3}^{\prime}u$, then it is a tree connecting $S$ in $G_{n}$. By Lemma 13, $\kappa(G_{n-1}^{3}\bigoplus G_{n-1}^{4})\geq r+2(n-2)\geq 4$. Hence, $G_{n}[(V(G_{n-1}^{3})\bigcup V(G_{n-1}^{4})\setminus \{v_{3}^{\prime}\}]$ is connected and it contains a tree $T_{r+2n-3}^{\prime}$ connecting $v_{1}^{\prime\prime}, v_{2}^{\prime\prime}$ and $v_{3}^{\prime\prime}$. Let $T_{r+2n-3}= T_{r+2n-3}^{\prime}\bigcup v_{1}v_{1}^{\prime\prime}\bigcup v_{2}v_{2}^{\prime\prime}\bigcup v_{3}v_{3}^{\prime\prime}$, then it is a tree connecting $S$ and the result holds. Otherwise, $u\in V(G_{n-1}^{4})$. Let $x$ be an in-neighbour of $v_{3}^{\prime}$ in $G_{n-1}^{3}$ such that one of the outside neighbour of $x$, say $z$, does not belong to $G_{n-1}^{4}$. This can be done as $r+2(n-2)\geq 4$. Hence, $G_{n}[\bigcup _{i\in[p_{n}]\setminus\{1,3,4\}}V(G_{n-1}^{i})]$ contains a tree, say $T_{r+2n-4}^{\prime}$ connecting $v_{1}^{\prime}, v_{2}^{\prime}$ and $z$. Let $T_{r+2n-4}=T_{r+2n-4}^{\prime}\bigcup v_{1}v_{1}^{\prime}\bigcup v_{2}v_{2}^{\prime}\bigcup zx\bigcup xv_{3}^{\prime}\bigcup v_{3}v_{3}^{\prime}$, then it is a tree connecting $S$ in $G_{n}$. By Lemma $13$, $G_{n}[(V(G_{n-1}^{3})\bigcup V(G_{n-1}^{4})\setminus \{v_{3}^{\prime}, x\}]$ is connected. Then there is a tree, say $T_{r+2n-3}^{\prime}$, connecting $v_{1}^{\prime\prime}, v_{2}^{\prime\prime}$ and $v_{3}^{\prime\prime}$. Let $T_{r+2n-3}=T_{r+2n-3}^{\prime}\bigcup v_{1}v_{1}^{\prime\prime}\bigcup v_{2}v_{2}^{\prime\prime}\bigcup v_{3}v_{3}^{\prime\prime}$, then it is a tree connecting $S$ in $G_{n}$. Combining the $T_{i}$s for $1\leq i \leq r+2n-3$, $r+2n-3$ internally disjoint trees connecting $S$ in $G_{n}$ are obtained.

Case $3$. Each copy contains at most one vertex of $M$.

Without loss of generality, suppose that $G_{n-1}^{2}, G_{n-1}^{3}, G_{n-1}^{4}$ contains $v_{1}^{\prime}, v_{2}^{\prime},  v_{3}^{\prime}$, respectively and  $G_{n-1}^{5}, G_{n-1}^{6}, G_{n-1}^{7}$ contains $v_{1}^{\prime\prime}, v_{2}^{\prime\prime},  v_{3}^{\prime\prime}$, respectively. As $G_{n}[\bigcup _{i=2}^{4}V(G_{n-1}^{i})]$ and $G_{n}[\bigcup _{i=5}^{7}V(G_{n-1}^{i})]$ as induced subgraphs of $G_{n}$ are both connected, there is a tree, say $T_{r+2n-4}^{\prime}$, connecting $v_{1}^{\prime}, v_{2}^{\prime}$ and $v_{3}^{\prime}$ in $G_{n}[\bigcup _{i=2}^{4}V(G_{n-1}^{i})]$ and a tree, say $T_{r+2n-3}^{\prime}$, connecting $v_{1}^{\prime\prime}, v_{2}^{\prime\prime}$ and $v_{3}^{\prime\prime}$ in $G_{n}[\bigcup _{i=5}^{7}V(G_{n-1}^{i})]$, respectively. Let $T_{r+2n-4}=T_{r+2n-4}^{\prime}\bigcup v_{1}v_{1}^{\prime}\bigcup v_{2}v_{2}^{\prime}\bigcup v_{3}v_{3}^{\prime}$ and $T_{r+2n-3}=T_{r+2n-3}^{\prime}\bigcup v_{1}v_{1}^{\prime\prime}\bigcup v_{2}v_{2}^{\prime\prime}\bigcup v_{3}v_{3}^{\prime\prime}$. Combining the $T_{i}$s for $1\leq i \leq r+2n-3$, $r+2n-3$ internally disjoint trees connecting $S$ in $G_{n}$ are obtained.
\hfill\qed

\begin{theorem}\label{thm2}
Let $G_{n}$ and $r$ be the same as definition $1$ and let  $G_{n}=G_{n-1}^{1}\bigoplus G_{n-1}^{2}\bigoplus \ldots\bigoplus \\G_{n-1}^{p_{n}}$. If any two vertices in different copies of $G_{n-1}$ have at most one common outside neighbour, then $\kappa_{3}(G_{n})=r+2n-3$, where $\kappa_{3}(G_{1})=r-1$.
\end{theorem}

\f {\bf Proof.} By definition $1$, $G_{n}$ is $[r+2(n-1)]$-regular. By Lemma~\ref{lem1}, $\kappa_{3}(G_{n})\leq \delta-1=r+2n-3$. To prove the result, we just need to show that $\kappa_{3}(G_{n})\geq r+2n-3$. We prove the result by induction on $n$.

Note that $\kappa_{3}(G_{1})=r-1$. Thus, the result holds for $n=1$. Next, assume that $n\geq 2$. Let $G_{n}=G_{n-1}^{1}\bigoplus G_{n-1}^{2} \bigoplus \ldots \bigoplus G_{n-1}^{p_{n}}$ and $v_{1}, v_{2}, v_{3}$ be any three distinct vertices of $G_{n}$. For convenience, let $S=\{v_{1}, v_{2}, v_{3}\}$ and the following three cases are considered.

Case $1$. $v_{1}, v_{2}$ and $v_{3}$ belong to the same copy of $G_{n-1}$.

Without loss of generality, let $S\subseteq V(G_{n-1}^{1})$. By the inductive hypothesis, there are $r+2n-5$ internally disjoint trees connecting $S$ in $G_{n-1}^{1}$. By Lemma~\ref{lem8}, there are $r+2n-3$ internally disjoint trees connecting $S$ in $G_{n}$ and the result is desired.

Case 2. $v_{1}, v_{2}$ and $v_{3}$ belong to two different copies of $G_{n-1}$.

Without loss of generality, let $v_{1}, v_{2}\in V(G_{n-1}^{1})$ and $v_{3}\in V(G_{n-1}^{2})$. By definition $1(6)$, $\kappa(G_{n-1}^{1})=r+2(n-2)$. Then there exist $r+2(n-2)$ internally disjoint paths $P_{1}, P_{2}, \ldots, P_{r+2(n-2)}$ between $v_{1}$ and $v_{2}$ in $G_{n-1}^{1}$. Let $H=G_{n-1}^{2} \bigoplus G_{n-1}^{3}\bigoplus\cdots \bigoplus G_{n-1}^{p_{n}}$. Then at most one outside neighbour of $v_{3}$ belongs to $V(G_{n-1}^{1})$ and the following two subcases are considered.

Subcase $2.1$. Neither of the two outside neighbours of $v_{3}$ belong to $G_{n-1}^{1}$, that is $d_{H}(v_{3})=r+2(n-1)$.

Choose $r+2(n-2)$ distinct vertices $x_{1}, x_{2}, \cdots, x_{r+2(n-2)}$ from $P_{1}, P_{2}, \ldots, P_{r+2(n-2)}$ such that $x_{i}\in V(P_{i})$ for $1\leq i\leq r+2(n-2)$, see Figure $5$. At most one of the paths has length $1$. If so, say $P_{1}$ and let $x_{1}=v_{1}$. Let $Y=\{x_{1}, x_{2}, \cdots, x_{r+2(n-2)}\} \bigcup \{v_{1}, v_{2}\}$. If $x_{1}\neq v_{1}$, let $Y^{\prime}=\{x^{\prime}|x^{\prime}$ is an outside neighbour of $x$ and $x\in Y\}$.  If $x_{1}= v_{1}$, let $Y^{\prime}=\{x^{\prime}|x^{\prime}$ is an outside neighbour of $x$ and $x\in Y\} \bigcup \{v_{1}^{\prime\prime}\}$, where $v_{1}^{\prime}$ and $v_{1}^{\prime\prime}$ are two outside neighbours of $v_{1}$. Clearly, $|Y|\geq r+2n-3$ and $|Y^{\prime}|=r+2(n-1)$. We can make sure that $|Y^{\prime}\bigcap G_{n-1}^{j}|\leq r+2(n-2)$ for each $j\in\{2, 3,\cdots, p_{n}\}$. If not, we can replace with the other outside neighbour of $x$ for some $x\in Y$. As $d_{H}(v_{3})=r+2(n-1)$. By Lemma~\ref{lem4}, there exist $r+2(n-1)$ internally disjoint $(v_{3}, Y^{\prime})$-paths $Q_{1}, Q_{2}, \cdots, Q_{r+2(n-1)}$ in $H$ such that the terminal vertex of $Q_{i}$ is $x_{i}^{\prime}$ for each $i\in \{1, 2,3,\cdots, r+2(n-2)\}$, the terminal vertex of $Q_{r+2n-3}$ is $v_{1}^{\prime}$ or $v_{1}^{\prime\prime}$ and the terminal vertex of $Q_{r+2n-2}$ is $v_{2}^{\prime}$. Let $T_{i}=P_{i}\bigcup Q_{i}\bigcup x_{i}x_{i}^{\prime}$ for $1\leq i \leq r+2(n-2), T_{r+2n-3}=Q_{r+2n-3}\bigcup Q_{r+2n-2}\bigcup v_{2}v_{2}^{\prime}\bigcup v_{1}v_{1}^{\prime}$ or $T_{r+2n-3}=Q_{r+2n-3}\bigcup Q_{r+2n-2}\bigcup v_{2}v_{2}^{\prime}\bigcup v_{1}v_{1}^{\prime\prime}$, then $r+2n-3$ internally disjoint trees connecting $S$ in $G_{n}$ are obtained.
\begin{figure}[!ht]
\begin{center}
\vskip1cm
\includegraphics[scale=0.8]{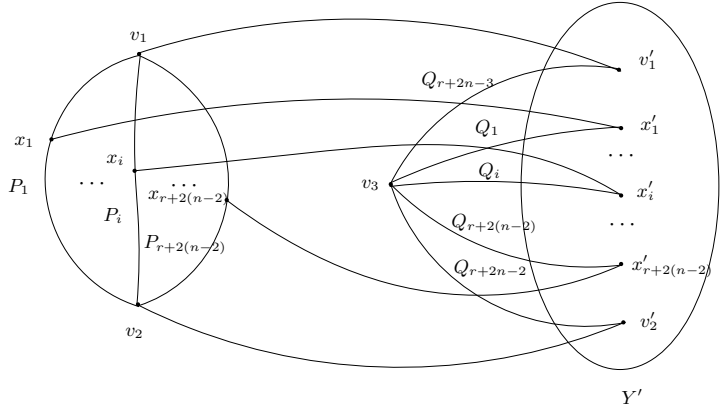}\label{F5}
\end{center}
\vskip0.1cm
\caption{Illustration of Subcase $2.1$}\label{F3}
\end{figure}

Subcase $2.2$. One of the outside neighbours of $v_{3}$ belongs to $G_{n-1}^{1}$, that is $d_{H}(v_{3})=r+2n-3$.

Without loss of generality, let $v_{3}^{\prime}$ be one of the outside neighbours of $v_{3}$ and belong to $G_{n-1}^{1}$. In addition, let $V(P)=\bigcup_{i=1}^{r+2(n-2)}V(P_{i})$.

If $v_{3}^{\prime}\notin V(P)$, as $G_{n-1}^{1}$ is connected, there is a $(v_{3}^{\prime}, v_{1})$-path $\widetilde{P}$ in $G_{n-1}^{1}$. Let $t$ be the first vertex of $\widetilde{P}$ which is in $V(P)$ and assume that $t\in V(P_{r+2(n-2)})$. Clearly, $P_{r+2(n-2)}\bigcup \widetilde{P}[v_{3}^{\prime},t]\\ \bigcup v_{3}v_{3}^{\prime}$ is a tree connecting $S$, denoted by $T_{r+2n-3}$. If $v_{3}^{\prime}\in V(P)$, without loss of generality, let $v_{3}^{\prime}\in V(P_{r+2(n-2)})$. Let $T_{r+2n-3}=P_{r+2(n-2)}\bigcup v_{3}v_{3}^{\prime}$, then it is a tree connecting $S$.

Next, choose $r+2n-5$ distinct vertices $x_{1}, x_{2}, \cdots, x_{r+2n-5}$ from $P_{1}, P_{2}, \ldots, P_{r+2n-5}$ such that $x_{i}\in V(P_{i})$ for $1\leq i\leq r+2n-5$. Denote $Y$ and $Y^{\prime}$ similarly as Subcase $2.1$. By Lemma~\ref{lem7} and the fact that $d_{H}(v_{3})=r+2n-3$, there exist $r+2n-3$ internally disjoint $(v_{3}, Y^{\prime})$-paths $Q_{1}, Q_{2}, \cdots, Q_{r+2n-3}$ in $H$ such that the terminal vertex of $Q_{i}$ is $x_{i}^{\prime}$ for each $i\in \{1, 2,3,\cdots, r+2n-5\}$, the terminal vertex of $Q_{r+2n-4}$ is $v_{1}^{\prime}$ or $v_{1}^{\prime\prime}$ and the terminal vertex of $Q_{r+2n-3}$ is $v_{2}^{\prime}$. Now, let $T_{i}=P_{i}\bigcup Q_{i}\bigcup x_{i}x_{i}^{\prime}$ for each $i\in\{1,2,\cdots r+2n-5\}, T_{r+2n-4}=Q_{r+2n-4}\bigcup Q_{r+2n-3}\bigcup v_{2}v_{2}^{\prime}\bigcup v_{1}v_{1}^{\prime}$ or $T_{r+2n-4}=Q_{r+2n-4}\bigcup Q_{r+2n-3}\bigcup v_{2}v_{2}^{\prime}\bigcup v_{1}v_{1}^{\prime\prime}$, then $r+2n-3$ internally disjoint trees connecting $S$ in $G_{n}$ are obtained.

Case 3. $v_{1}, v_{2}$ and $v_{3}$ belong to three different copies of $G_{n-1}$, respectively.

Without loss of generality, we assume that $v_{1} \in V(G_{n-1}^{1}), v_{2}\in V(G_{n-1}^{2})$ and $v_{3}\in V(G_{n-1}^{3})$. Let $W=\{v_{1}^{\prime}, v_{2}^{\prime}, v_{3}^{\prime}, v_{1}^{\prime\prime}, v_{2}^{\prime\prime}, v_{3}^{\prime\prime}\}$, where $v_{i}^{\prime}$ and $v_{i}^{\prime\prime}$ are the two outside neighbours of $v_{i}$ for $1\leq i \leq 3$. The following three subcases are considered.

Subcase $3.1$. $W\subseteq V(G_{n-1}^{1})\bigcup V(G_{n-1}^{2})\bigcup V(G_{n-1}^{3})$.

Let $H=G_{n-1}^{1}\bigoplus G_{n-1}^{2}$. Since one of the two outside neighbours of $v_{1}$ belongs to $G_{n-1}^{2}$ and one of the two outside neighbours of $v_{2}$ belongs to $G_{n-1}^{1}$. Hence, $d_{H}(v_{1})=d_{H}(v_{2})=r+2n-3$. By Lemma~\ref{lem7}, there exist $r+2n-3$ internally disjoint paths $P_{1}, P_{2}, \ldots, P_{r+2n-3}$ between $v_{1}$ and $v_{2}$ in $H$. Let $v_{3}^{\prime}$ be an outside neighbour of $v_{3}$, then  $v_{3}^{\prime}\in V(H)$. Let $V(P)=\bigcup_{i=1}^{r+2n-3}V(P_{i})$, as $H$ is connected, there is a path $\widetilde{P}$ from $v_{3}^{\prime}$ to $v_{1}$ in $H$. Let $t$ be the first vertex of $\widetilde{P}$ which is in $V(P)$ and assume that $t\in V(P_{r+2n-3})$. Clearly, $P_{r+2n-3}\bigcup \widetilde{P}[v_{3}^{\prime},t]\bigcup v_{3}v_{3}^{\prime}$ contains a tree connecting $S$, denoted by $T_{r+2n-3}$. If $v_{3}^{\prime}\in V(P)$, then let $v_{3}^{\prime}\in V(P_{r+2n-3})$ and $T_{r+2n-3}=P_{r+2n-3}\bigcup v_{3}v_{3}^{\prime}$, then it is a tree connecting $S$.

Let $x_{i}\in V(P_{i})\bigcap N_{H}(v_{1})$ for each $i\in\{1, 2,\cdots, r+2n-4\}$. If the outside neighbour of $v_{1}$ in $H$ does not belong to $x_{i}$s for $1\leq i\leq r+2n-4$,  let $X=\{x_{1}, x_{2},\cdots, x_{r+2n-4}\}$. If the outside neighbour of $v_{1}$ in $H$ belongs to $x_{i}$s for $1\leq i\leq r+2n-4$, say $x_{1}$, and let $X=\{v_{1}, x_{2},\cdots, x_{r+2n-4}\}$. Then $X\subseteq V(G_{n-1}^{1})$ and $|X|=r+2n-4$.
Let $H^{\prime}=G_{n-1}^{3}\bigoplus G_{n-1}^{4}\bigoplus \cdots \bigoplus G_{n-1}^{p_{n}}$ and $x_{i}^{\prime}$ be one of the two outside neighbours of $x_{i}$ such that $x_{i}^{\prime}\in V(H^{\prime})$ for each $i\in \{1,2,\cdots, r+2n-4\}$.

If $X=\{x_{1}, x_{2},\cdots, x_{r+2n-4}\}$, let $X^{\prime}=\{x_{1}^{\prime}, x_{2}^{\prime}, \cdots, x_{r+2n-4}^{\prime}\}$. By Lemma~\ref{lem1}, $|X^{\prime}|=r+2n-4$. As $d_{H^{\prime}}(v_{3})=r+2n-4$, by Lemma~\ref{lem4}, there exist $r+2n-4$ internally disjoint $(v_{3}, X^{\prime})$-paths $Q_{1}, Q_{2}, \cdots, Q_{r+2n-4}$ in $H^{\prime}$ such that the terminal vertex of $Q_{i}$ is $x_{i}^{\prime}$ for each $i\in \{1, 2,3,\cdots, r+2n-4\}$. Note that at most one of $Q_{i}$s for $1\leq i\leq r+2n-4$ has length one.
Let $T_{i}=P_{i}\bigcup Q_{i}\bigcup x_{i}x_{i}^{\prime}$ for $1\leq i \leq r+2n-4$. Combining with $T_{i}$s for $1\leq i\leq r+2n-3$, then $r+2n-3$ internally disjoint trees connecting $S$ in $G_{n}$ are obtained.

If $X=\{v_{1}, x_{2},\cdots, x_{r+2n-4}\}$, let $X^{\prime}=\{v_{1}^{\prime}, x_{2}^{\prime}, \cdots, x_{r+2n-4}^{\prime}\}$, where $v_{1}^{\prime}\in V(H^{\prime})$. With the similar method as $X=\{x_{1}, x_{2},\cdots, x_{r+2n-4}\}$, $r+2n-3$ internally disjoint trees $T_{i}$s for $1\leq i\leq r+2n-3$ connecting $S$ in $G_{n}$ can be obtained.

Subcase $3.2$. $W \nsubseteq V(G_{n-1}^{1})\bigcup V(G_{n-1}^{2})\bigcup V(G_{n-1}^{3})$.

Since $W \nsubseteq V(G_{n-1}^{1})\bigcup V(G_{n-1}^{2})\bigcup V(G_{n-1}^{3})$, at least one of the outside neighbours of $v_{3}$ does not belong to $V(G_{n-1}^{1})\bigcup V(G_{n-1}^{2})$. Let $H=G_{n-1}^{1}\bigoplus G_{n-1}^{2}$ and $H^{\prime}=G_{n-1}^{3}\bigoplus G_{n-1}^{4}\\ \bigoplus \cdots \bigoplus G_{n-1}^{p_{n}}$. Then select $r+2n-4$ vertices from $G_{n-1}^{1}\setminus \{v_{1}\}$, say $x_{1}, x_{2}, \cdots, x_{r+2n-4}$, such that one of the outside neighbours $x_{i}^{\prime}$ of $x_{i}$ belongs to $G_{n-1}^{2}$ for each $i\in\{1,2,\cdots, r+2n-4\}$. Further, we request that $x_{i}$ and $v_{2}$ have different outside neighbours for $1\leq i \leq r+2n-4$.

Let $S=\{x_{1}, x_{2}, \cdots, x_{r+2n-4}\}$ and $S^{\prime}=\{x_{1}^{\prime}, x_{2}^{\prime}, \cdots, x_{r+2n-4}^{\prime}\}$. By definition $1(6)$, $\kappa(G_{n-1}^{1})=\kappa(G_{n-1}^{2})=r+2n-4$. By Lemma~\ref{lem4}, there exist $r+2n-4$ internally disjoint $(v_{1}, S)$-paths $P_{1}, P_{2}, \ldots, P_{r+2n-4}$ in $G_{n-1}^{1}$ such that the terminal vertex of $P_{i}$ is $x_{i}$ and there exist $r+2n-4$ internally disjoint $(v_{2}, S^{\prime})$-paths $P_{1}^{\prime}, P_{2}^{\prime}, \ldots, P_{r+2n-4}^{\prime}$ in $G_{n-1}^{2}$ such that the terminal vertex of $P_{i}^{\prime}$ is $x_{i}^{\prime}$ for $1\leq i\leq r+2n-4$. Thus, we obtain $r+2n-4$ internally disjoint paths between $v_{1}$ and $v_{2}$ in $H$, where $\widetilde{P_{i}}=v_{1}P_{i}x_{i}x_{i}^{\prime}P_{i}^{\prime}v_{2}$ for each $i\in\{1,2,\cdots, r+2n-4\}$.

Now, let $v_{i}^{\prime\prime}$ be one of the outside neighbours of $v_{i}$ such that $v_{i}^{\prime\prime}\in V(H^{\prime})$ for $i=1,2$ and $x_{i}^{\prime\prime}$ be the other outside neighbour of $x_{i}$ such that $x_{i}^{\prime\prime}\in V(H^{\prime})$ for $1\leq i\leq r+2n-4$. Let
$Y=\{x_{1}^{\prime\prime}, x_{2}^{\prime\prime}, x_{3}^{\prime\prime}, \cdots, x_{r+2n-4}^{\prime\prime}, v_{1}^{\prime\prime}, v_{2}^{\prime\prime}\}$. Then $Y\subseteq V(H^{\prime})$ and $|Y|\geq r+2n-3$. If $v_{1}^{\prime\prime}\neq v_{2}^{\prime\prime}$, then $|Y|=r+2n-2$. If $v_{1}^{\prime\prime}=v_{2}^{\prime\prime}$, then $|Y|=r+2n-3$.

Subcase $3.2.1$. Neither of the two outside neighbours of $v_{3}$ belong to $\bigcup _{i=1}^{2}V(G_{n-1}^{i})$.

In this case, $d_{H^{\prime}}(v_{3})=r+2n-2$. If $|Y|=r+2n-2$, the proof is similar to Subcase $2.1$. If $|Y|=r+2n-3$, the proof is also similar to Subcase $2.1$ except that the paths $ Q_{r+2n-3}$ and $Q_{r+2n-2}$ become the same path.

Subcase $3.2.2$. One of the two outside neighbours of $v_{3}$ belongs to $\bigcup _{i=1}^{2}V(G_{n-1}^{i})$.

In this case, $d_{H^{\prime}}(v_{3})=r+2n-3$. If $|Y|=r+2n-2$, the proof is similar to Subcase $2.2$. If $|Y|=r+2n-3$, the proof is also similar to Subcase $2.2$ except that the paths $ Q_{r+2n-4}$ and $Q_{r+2n-3}$ become the same path.

Hence, $r+2n-3$ internally disjoint trees connecting $S$ in $G_{n}$ can be obtained in any case.
\hfill\qed
\section{Applications to some interconnection networks}
\subsection{Application to the alternating group graph $AG_{n}$}
\begin{lem}{\rm (\cite{ZHAO})}\label{lem9}
Let $AG_{n}=AG_{n-1}^{1}\bigoplus AG_{n-1}^{2}\bigoplus \ldots\bigoplus AG_{n-1}^{n}$ for $n\geq 3$. Then any two vertices in different copies of $AG_{n-1}$ have at most one common outside neighbour.
\end{lem}

\begin{cor}\label{cor4.5}
$\kappa_{3}(AG_{n})=2n-5$ for $n\geq 3$.
\end{cor}

\f {\bf Proof.} By definition 1,$AG_{n}$ can be regarded as the special regular graph $G_{n-2}$ with $G_{1}=AG_{3}$, $a=3$, $r=2$, $s=2$, $p_{n-2}=n$ and $N=ap_{2}p_{3}\cdots p_{n-2}=\frac{n!}{2}$. By Lemma 2, $\kappa(AG_{3})=2$. By Lemma~\ref{lem1}, $\kappa_{3}(AG_{3})\leq 1$. By Lemma~\ref{lem2}, $\kappa_{3}(AG_{3})\geq 1$. Thus, $\kappa_{3}(AG_{3})=1$. Thus, by Lemma $17$ and Theorem $1$, $\kappa_{3}(AG_{n})=2n-5$ for $n\geq 3$.

\subsection{Application to the $k$-ary $n$-cube $Q_{n}^{k}$}
\begin{lem}\label{lem9.19}
Let $Q_{n}^{k}=Q_{n-1}^{k}[0]\bigoplus Q_{n-1}^{k}[1]\bigoplus \ldots\bigoplus Q_{n-1}^{k}[k-1]$ for $k\geq 3$ and $n\geq 1$. Then any two vertices in different copies of $Q_{n-1}^{k}$ have at most one common outside neighbour.
\end{lem}

\f {\bf Proof.} Let $u, v\in V(Q_{n}^{k}), u\neq v$ and they belong to different copies of $Q_{n-1}^{k}$. Without loss of generality, let $u=u_{1}u_{2}u_{3}\cdots u_{n-1}0\in V(Q_{n-1}^{k}[0])$ and $v=v_{1}v_{2}v_{3}\cdots v_{n-1}1\in V(Q_{n-1}^{k}[1])$. Then the two outside neighbours of $u$ are $u^{\prime}=u_{1}u_{2}u_{3}\cdots u_{n-1}1$ and $u^{\prime\prime}=u_{1}u_{2}u_{3}\cdots u_{n-1}(k-1)$, and the two outside neighbours of $v$ are $v^{\prime}=v_{1}v_{2}v_{3}\cdots v_{n-1}0$ and $v^{\prime\prime}=v_{1}v_{2}v_{3}\cdots v_{n-1}2$. If $u$ and $v$ have two common outside neighbours, then $\{u^{\prime}, u^{\prime\prime}\}=\{v^{\prime}, v^{\prime\prime}\}$. As $u^{\prime}\neq v^{\prime}$, then $u^{\prime}=v^{\prime\prime}$ and $v^{\prime}=u^{\prime\prime}$. However, $u^{\prime}\neq v^{\prime\prime}$ clearly, which is a contradiction. Thus, $u$ and $v$ have at most one common outside neighbour.
\hfill\qed

\begin{cor}\label{cor4.6}
$\kappa_{3}(Q_{n}^{k})=2n-1$ for $k\geq 3$ and $n\geq 1$.
\end{cor}
\f {\bf Proof.} By definition 1, $Q_{n}^{k}(k\geq 3)$ can be regarded as the special regular graph $G_{n}$ with $G_{1}=Q_{1}^{k}$, $a=k$, $r=2$, $s=2$, $p_{n}=k$ and $N=ap_{2}p_{3}\cdots p_{n}=k^{n}$. By Lemma 4, $\kappa(Q_{1}^{k})=2$. By Lemma~\ref{lem1}, $\kappa_{3}(Q_{1}^{k})\leq 1$. By Lemma~\ref{lem2}, $\kappa_{3}(Q_{1}^{k})\geq 1$. Thus, $\kappa_{3}(Q_{1}^{k})=1$. By Lemma~\ref{lem9.19} and Theorem $1$, $\kappa_{3}(Q_{n}^{k})=2n-1$ for $k\geq 3$ and $n\geq 1$.

\subsection{Application to the split-star network $S_{n}^{2}$}
\begin{lem}\label{lem9.1}
Let $S_{n}^{2}=S_{n-1}^{2}[1]\bigoplus S_{n-1}^{2}[2]\bigoplus \ldots\bigoplus S_{n-1}^{2}[n]$ for $n\geq 3$. Then any two vertices in different copies of $S_{n-1}^{2}$ have at most one common outside neighbour.
\end{lem}

\f {\bf Proof.} Let $u, v\in V(S_{n}^{2}), u\neq v$ and they belong to different copies of $S_{n-1}^{2}$. Without loss of generality, let $u=u_{1}u_{2}u_{3}\cdots u_{n-1}1\in V(S_{n-1}^{2}[1])$ and $v=v_{1}v_{2}v_{3}\cdots v_{n-1}2\in V(S_{n-1}^{2}[2])$. Then the two outside neighbours of $u$ are $u^{\prime}=u(1n)=1u_{2}u_{3}\cdots u_{n-1}u_{1}$ and $u^{\prime\prime}=u(2n)=u_{1}1u_{3}\cdots u_{n-1}u_{2}$, and the two outside neighbours of $v$ are $v^{\prime}=v(1n)=2v_{2}v_{3}\cdots v_{n-1}v_{1}$ and $v^{\prime\prime}=v(2n)=v_{1}2v_{3}\cdots v_{n-1}v_{2}$. If $u$ and $v$ have two common outside neighbours, then $\{u^{\prime}, u^{\prime\prime}\}=\{v^{\prime}, v^{\prime\prime}\}$. As $u^{\prime}\neq v^{\prime}$, then $u^{\prime}=v^{\prime\prime}$ and $v^{\prime}=u^{\prime\prime}$. By $u^{\prime}=v^{\prime\prime}$, we have that $v_{1}=1$ and $u_{2}=2$. By $u^{\prime\prime}=v^{\prime}$, we have that $v_{2}=1$ and $u_{1}=2$. That is, $u_{1}=u_{2}=2$, which is a contradiction. Thus, $u$ and $v$ have at most one common outside neighbour.
\hfill\qed

\begin{cor}\label{cor4.7}
$\kappa_{3}(S_{n}^{2})=2n-4$ for $n\geq 3$.
\end{cor}
\f {\bf Proof.} By definition 1, $S_{n}^{2}$ can be regarded as the special regular graph $G_{n-2}$ with $G_{1}=S_{3}^{2}$, $a=6$, $r=3$, $s=2$, $p_{n-2}=n$ and $N=ap_{2}p_{3}\cdots p_{n-2}=n!$. By Lemma 6, $\kappa(S_{3}^{2})=3$. By Lemma~\ref{lem1}, $\kappa_{3}(S_{3}^{2})\leq 2$. By Lemma~\ref{lem2}, $\kappa_{3}(S_{3}^{2})\geq 2$. Thus, $\kappa_{3}(S_{3}^{2})=2$. Thus, by Lemma~\ref{lem9.1} and Theorem $1$, $\kappa_{3}(S_{n}^{2})=2n-4$ for $n\geq 3$.

\subsection{Application to the bubble-sort-star network $BS_{n}$}
\begin{lem}\label{lem9}
Let $BS_{n}=BS_{n-1}^{1}\bigoplus BS_{n-1}^{2} \bigoplus \ldots \bigoplus BS_{n-1}^{n}$. Then any two vertices in different copies of $BS_{n-1}$ have at most one common outside neighbour.
\end{lem}
\f {\bf Proof.} Let $u, v\in V(BS_{n}), u\neq v$ and they belong to different copies of $BS_{n-1}$. Without loss of generality, let $u=u_{1}u_{2}u_{3}\cdots u_{n-1}1\in V(BS_{n-1}^{1})$ and $v=v_{1}v_{2}v_{3}\cdots v_{n-1}2\in V(BS_{n-1}^{2})$. Then the two outside neighbours of $u$ are $u^{\prime}=u(1n)=1u_{2}u_{3}\cdots u_{n-1}u_{1}$ and $u^{\prime\prime}=u(n-1,n)=u_{1}u_{2}\cdots u_{n-2}1u_{n-1}$, and the two outside neighbours of $v$ are $v^{\prime}=v(1n)=2v_{2}v_{3}\cdots v_{n-1}v_{1}$ and $v^{\prime\prime}=v(n-1,n)=v_{1}v_{2}\cdots v_{n-2}2v_{n-1}$. If $u$ and $v$ have two common outside neighbours, then $\{u^{\prime}, u^{\prime\prime}\}=\{v^{\prime}, v^{\prime\prime}\}$. As $u^{\prime}\neq v^{\prime}$, then $u^{\prime}=v^{\prime\prime}$ and $v^{\prime}=u^{\prime\prime}$. By $u^{\prime}=v^{\prime\prime}$, we have that $v_{1}=1$ and $u_{n-1}=2$. By $u^{\prime\prime}=v^{\prime}$, we have that $v_{n-1}=1$ and $u_{1}=2$. That is, $u_{1}=u_{n-1}=2$, which is a contradiction. Thus, $u$ and $v$ have at most one common outside neighbour.
\hfill\qed

\begin{cor}\label{cor4.8}
$\kappa_{3}(BS_{n})=2n-4$ for $n\geq 3$.
\end{cor}

\f {\bf Proof.} By definition 1, $BS_{n}$ can be regarded as the special regular graph $G_{n-2}$ with $G_{1}=BS_{3}$, $a=6$, $r=3$, $s=2$, $p_{n-2}=n$ and $N=ap_{2}p_{3}\cdots p_{n-2}=n!$. By Lemma 8, $\kappa(BS_{3})=3$. By Lemma~\ref{lem1}, $\kappa_{3}(BS_{3})\leq 2$. By Lemma~\ref{lem2}, $\kappa_{3}(BS_{3})\geq 2$. Thus, $\kappa_{3}(BS_{3})=2$. Thus, by Lemma~\ref{lem9.1} and Theorem $1$, $\kappa_{3}(BS_{n})=2n-4$ for $n\geq 3$.
\hfill\qed

\section{Concluding remarks}
The generalized $k$-connectivity is a generalization of traditional connectivity. In this paper, we focus on the recursively constructed graphs $G_{n}$ with two outside neighbours. As applications of the main result, the generalized $3$-connectivity of many famous networks such as the alternating group graph $AG_{n}$, the $k$-ary $n$-cube $Q_{n}^{k}$, the split-star network $S_{n}^{2}$ and the bubble-sort-star graph $BS_{n}$ etc. can be obtained directly. The corresponding results for the recursively constructed graphs with any outside neighbours will be interesting and challenged problems which can be studied in the future.

\section*{Acknowledgments}
This work was supported by the National Natural Science Foundation of China (No.11731002), the Fundamental Research Funds for the Central Universities (No. 2016JBM071, 2016JBZ012) and the $111$ Project of China (B16002). The authors would like to thank the editor and the anonymous reviewers for their kind suggestions which greatly improved the original manuscript.

\end{document}